# Optimizing Electric Carsharing System Operations and Battery Management: Integrating V2G, B2G and Battery Swapping Strategies

Shuang Yang[a], Gonçalo Homem de Almeida Correia[b], Jianjun Wu[a], Huijun Sun[a,*]

[a] School of Systems Science, Beijing Jiaotong University, Beijing, China

[b] Department of Transport & Planning, Delft University of Technology, Delft, the Netherlands

**Abstract**

Shared electric vehicles (SEVs) have emerged as a promising solution to contribute to sustainable urban mobility. However, ensuring the efficient operation and effective battery management of SEV systems remains a complex challenge. This challenge stems from factors such as slow plug-in charging, the potential role of SEVs in balancing grid load pressure, and the optimization of SEV operations to ensure their economic viability. To tackle these challenges, this paper introduces an integrated strategy for optimizing various aspects of SEV systems, encompassing strategies like Vehicle-to-Grid (V2G), Battery-to-Grid (B2G), and battery swapping. This approach is built on a space-time-energy network model that facilitates the optimization of battery charging and discharging scheduling, SEV operations like relocations and battery swapping, battery swapping station selection and the number of batteries. The objective of this approach is to maximize profits while addressing operational constraints and the complexities of energy management within SEV systems. Given the substantial complexity that arises with large-problem scales, the paper introduces a column generation-based heuristic algorithm. Extensive experimental validation is conducted, including sensitivity analysis on different charging speeds and fleet sizes. The results illuminate the impact of varying charging rates and fleet sizes on performance indicators. Notably, it is observed that battery swapping is particularly effective as an auxiliary charging method when the number of vehicles is limited. Conversely, in scenarios with a large fleet, the necessity for battery swapping diminishes. Moreover, results show the effectiveness of V2G and B2G technologies in grid load balancing. Additionally, it is found that when opting for a recharging facility with slow charging speed, utilizing battery swapping proves to be a valuable method to compensate

---

[*] Corresponding author.

*E-mail address*: hjsun1@bjtu.edu.cn (H.J. Sun).

for the shortage of available charged vehicles.

## 1. Introduction

Carsharing services, as a form of car rental, offer citizens the opportunity to use cars without the burden of ownership. The entire process, spanning from reserving the vehicle to its collection, return and subsequent billing, is typically handled through a user-friendly mobile app. This concept not only liberates users from the finance encumbrances linked with car ownership, such as insurance premiums and parking expenditures, but also provides flexibility and convenience, contributing to its rapid growth (Mounce and Nelson, 2019). As a result, carsharing holds great promise as one of the preferred modes of travel in the future.

In recent years, the urgency to reduce greenhouse gas emissions has intensified due to the worsening environmental crisis. Statistics show that electric vehicles (EVs) produce greenhouse gas emissions approximately 17%-30% lower than petrol and diesel cars, proving that EVs help reduce negative impacts on climate change and air quality (European Environment Agency, 2018). Hence, governmental bodies across the globe have taken proactive measures by enacting various supportive policies tailored to bolster the adoption of EVs.

The use of electric mobility and carsharing has given birth to a symbiotic fusion of them. Shared electric vehicles (SEVs) integrate the convenience and accessibility intrinsic to carsharing with the commendable ecological merits inherent in EVs. The combination represents a significant step toward creating a more sustainable and efficient transportation system. Many carsharing companies have already channeled investments into the operation of SEVs, like EVCARD in China and Zity in Europe. However, the growth of SEVs faces challenges arising from the inherent limitations of EVs. Range anxiety and long charging duration, particularly when using slow and cheaper charging solutions, are major issues in the operation of SEVs.

To improve the operation efficiency, operators must optimize the energy replenishment process to ensure that EVs are charged with electricity timely. Currently, there are two primary methods of recharging EVs: battery swapping and plug-in charging. Among these options, plug-in charging takes precedence because of its cost-effectiveness and adaptability.

However, it is not devoid of limitations, notably the long charging time, even when utilizing fast charging technologies. Conversely, battery swapping offers the distinct advantage of quickly replacing depleted batteries with fully charged batteries for vehicles, making it suitable for long-distance travel or emergency situations. The advent of battery swapping has revolutionized the replenishment time for EVs, aligning it with the rapid refueling capabilities of conventional fossil-fuel cars. This is critical to time-sensitive services like carsharing. Moreover, carsharing platforms are innately well-suited to facilitate battery pack standardization, as they can enforce uniform battery pack requirements across their fleet. Currently, some shared mobility companies like Zbee and Gogoro, have embraced battery swapping services.

The difference between implementing and not implementing battery swapping in carsharing systems lies in two main aspects. Firstly, implementing battery swapping requires dedicated stations with specialized facilities for battery exchange. Another distinction arises in the number of batteries in the system. In the scenario where battery swapping is not employed, the number of batteries in the system matches the number of electric vehicles. If battery swapping is utilized, additional batteries need to be purchased and stored at the dedicated stations. It is worth noting that the cost of building dedicated stations and purchasing supplementary batteries constitutes a significant investment. This elevated cost factor poses a significant barrier to the widespread adoption of battery swapping. As a result, battery swapping is often perceived as a complementary approach to plug-in charging.

Regrettably, as we will elaborate, a gap exists in terms of comprehensive models that factor in both recharging techniques and their integration within the operations of SEVs. The modelling deficiency inherently encompasses the decision of whether to adopt the battery swapping approach, and if so, the determination of the required number of swappable batteries, their distribution, and activity schedules, taking into account both fixed and operational costs. Addressing this modelling gap would enable carsharing operators to refine their recharging operations and achieve a balance between costs and charging efficiency.

As dependence on the power grid continues to grow, the grid is under increasing strain, particularly during peak periods. To address the challenges posed by peak-load demand, governments have instituted various strategies, including the adoption of time-of-use pricing, resulting in higher prices during peak hours. Furthermore, in this context, EVs can also play a role in energy management through vehicle-to-grid (V2G) technology.

V2G is a bidirectional energy exchange system that enables EVs to feed stored energy back into the grid and draw energy from the grid when needed. During peak grid demands, EVs can discharge energy into the grid, easing power demand pressures and improving grid stability. An additional innovation is the direct bidirectional connection of batteries to the grid, termed Battery-to-Grid (B2G), which functions independently of EVs. V2G and B2G yield several benefits such as saving costs of EV owners, and enhancing the flexibility and stabilization of grid (Noel et al., 2019; Liao et al., 2021).

If SEVs are equipped with V2G and their batteries with B2G, SEVs and batteries have the potential not only to enhance grid stability but also generate additional revenue for operators. This integration brings forth a scenario where EVs and batteries simultaneously participate in both the carsharing service system and the power transmission system. Within this context, the complicated coupling of users, EVs and batteries presents a great operational challenge for operators. Mishandling this intricate interplay may result in reduced system efficiency, diminished user satisfaction, and a decline in profits. As a result, effectively coordinating the scheduling of vehicles and batteries, while considering this intricate interrelationship, becomes crucial.

Another well-known challenge in carsharing systems, especially within one-way systems, is the issue of vehicle imbalance (Lu et al., 2021; Huang et al., 2018). One-way systems provide users with the flexibility to pick up and return vehicles at different locations. However, due to fluctuations in demand, it is common to encounter situations where some stations confront vehicle shortages while others have an excess of vehicles. Such an imbalance negatively affects the vehicle utilization. In response to this challenge, vehicle relocation has emerged as a widely recognized and efficacious strategy. The process involves moving excess vehicles from stations with high inventories to those experiencing shortages, thereby achieving a more balanced distribution of vehicle. In order to maximize system efficiency, relocation should be jointly optimized with the recharging operation. Without synchronized optimization, simply relocating vehicles might inadvertently delay recharging and potentially affect the overall system efficiency. This integration can potentially improve the utilization efficiency of recharging infrastructure and increase the vehicle availability with sufficient battery power.

To address the above issues, this study proposes a space-time-energy network flow model. This model aims to optimize strategic decisions concerning battery swapping station siting, determining the number and distribution of swappable batteries, as well as operational decisions including charging, discharging, battery swapping of EVs and swappable batteries, and EV relocations. Notably, this optimization framework considers a range of charging methodologies, including battery swapping, V2G, and B2G. As far as the authors know, this research is the first attempt to tackle this integrated optimization problem. Recognizing the computational complexity of this model, we introduce a column generation-based (CG-based) heuristic algorithm. To show the effectiveness and validity of both the model and algorithm, we carry out numerical experiments using synthetic case studies and real network data, inclusive of sensitivity analysis.

This paper is organized as follows: Section 2 provides an overview of related work and highlights the contribution of this paper. This is followed by Section 3, where the formulation of the proposed optimization model based on the space-time-energy network flow model is introduced. Section 4 outlines our CG-based solution approach designed to solve the proposed model. To assess computation efficiency and gain insights into system operations, we conduct numerical studies in Section 5. Lastly, Section 6 provides the main conclusions from the study and outlines potential directions for future research.

## 2. Literature review

We review the related literature on the topics of SEVs and battery swapping. This review aims to not only provide an up-to-date and comprehensive understanding of the methodology for strategic and operational problems of SEVs and battery swapping but also shed light on areas where further investigation is needed, clearly positioning the scope and focus of the present study.

*2.1 Shared electric vehicles*

2.1.1 Strategic and operational problems of shared electric vehicles

In recent years, shared electric mobility has drawn growing attention. The literature on the topic of SEVs is extensive and can be categorized into strategic-level and operational-level planning problems (Ait-Ouahmed et al., 2018; Liu et al., 2022; Boyaci et al., 2015). Strategic-level decisions focus on long-time planning strategies, including the number, locations, capacity of parking and charging stations, fleet and staff size, as well as the initial vehicle distribution (Xu et al., 2018; Li et al., 2016; Zhao et al., 2021). For example, Li et al. (2016) propose a comprehensive design network for SEV systems, jointly optimizing the station location and EV fleet size under uncertain demand. A continuum approximation approach is introduced to solve large-scale problem instances. Another study by Zhao et al. (2021) also explores the deployment of charging stations and fleet operation issues of a station-based SEV system. They develop a simulation-based optimization model along with a tailored heuristic algorithm to tackle this model. Readers may refer to the review articles by Yao et al. (2022) for more details on strategic decision-making in SEV systems.

When it comes to the operational level, studies concentrate on the daily management and operational problems, particularly concerning relocation and charging operations. Often, these two issues are addressed simultaneously in the literature. An inherent concern within one-way station-based systems is vehicle imbalance, which refers to the imbalance in vehicle supply and demand between different regions or time periods (Nourinejad et al., 2015; Kek et al., 2009; Yang et al., 2021a). To tackle this challenge, operators usually perform relocations. Two relocation strategies have been proposed and implemented: operator-based and user-based relocations. Operator-based relocations involve hiring staff to move vehicles between stations (Nourinejad et al., 2015; Zhao et al., 2018; Yang et al., 2021a; Chang et al., 2022; Bruglieri et al., 2019). On the other hand, the user-based relocation strategy encourages users to actively participate in relocations by offering incentives such as discounts or rewards (Jorge et al., 2015; Wang et al., 2021; Zhan et al., 2022; Di Febbraro et al., 2019; Zhang et al., 2022; Stokkink and Geroliminis, 2021). Operator-based relocations are generally more straightforward and easier to implement compared to user-based approaches but they are also quite costly (Santos and Correia, 2019). This paper specifically focuses on the operator-based relocation strategy.

The utilization of EVs presents challenges in managing relocations due to the necessity of monitoring the state of charge (SOC) of the vehicles. Ensuring sufficient remaining battery power for upcoming rentals or relocations is vital to maintaining the operational efficiency of SEV systems. In the current literature, some simplifications are made regarding the recharging process when jointly optimizing recharging and relocation. This simplification assumes that once vehicles are returned to stations, they are immediately charged and must remain in this charged state for a predetermined duration referred to as "mandatory fixed charge". Alternatively, they remain charging until fully replenished before being designated for the subsequent trip or relocation, a concept termed "mandatory full charge" (Ait-Ouahmed et al., 2018; Boyaci et al., 2015; Xu et al., 2018; Xu et al., 2020; Brandstätter et al., 2017). Regrettably, this assumption restricts the utilization of partially charged vehicles especially during peak demand periods. Consequently, this practice could lead to the underutilization of vehicles, an overestimation of the required fleet size and a reduction in efficiency within carsharing systems. To address this issue, there arises a need to consider more adaptable and responsive charging strategies that can optimize energy consumption and enhance infrastructure efficiency. To fulfill this need, Gambella et al. (2018) and Zhao et al. (2018) use the time-space network to represent activity arcs involving user travel, staff movement, and EV travel within the context of integrated EV rebalancing and staff relocation. Jamshidi et al. (2021) propose a sequential mixed integer linear programming approach to determine the dynamic planning for recharging and relocations of electric taxis. To obtain these executable decisions for individual vehicles, they also adopt the time-space network to track the status of each vehicle. In their directed graph, the daily operational period is divided into a number of small intervals. Each node corresponds to a station during a specific time interval. They meticulously factor in EV battery consumption and the recharging process within their daily operational optimization model. They track the SOC of each vehicle by explicitly defining a time-dependent state variable. However, it is worth noting that the approach adopted by Gambella et al. (2018) and Zhao et al. (2018) involves a comprehensive depiction of the activities of each vehicle and staff member. Unfortunately, this level of detail results in complex, large-scale models that prove to be computationally intractable when applied to real-world carsharing systems.

Apart from the time-space network, a comprehensive representation method for vehicle status has gained attention in the literature, termed the space-time-energy network (Gonçalves Duarte Santos et al., 2023; Bekli et al., 2021; Zhang et al., 2019; Chen and Liu, 2022; Zhang et al., 2021). This approach does not involve describing each individual vehicle but aggregates vehicles with the same battery status. The space-time-energy network flow model serves as an extension of the conventional space-time network flow model. In addition to considering stations and time intervals, this modelling approach discretizes battery capacity into several levels, enabling the monitoring of SOC over time. Each node in the network is represented by a triple tuple: station, time and energy. Zhang et al. (2019) is the first to adopt the space-time-energy network flow model to formulate the EV assignment problem with the decision of EV routing, relocation and charging. Recognizing the uncertainties in demand and the existence of multiple charging outlets, Chen and Liu (2022) expand upon the work of Zhang et al. (2019) by proposing an event-activity space-time-energy network. This enhanced model explicitly captures diverse characteristics of different charging facilities and enables the tracking of charging choices and battery levels of vehicles.

Expanding the consideration to various types of charging outlets, Gonçalves Duarte Santos et al. (2023) also incorporate different types of vehicles with varying battery capacities and energy consumption rates. They develop the time-space-energy model to determine the movement and charging of EVs, as well as the required number of vehicles and charging facilities.

Using the space-time-energy network flow model avoids the need to individually model each vehicle, thus reducing computational complexity to some extent. Additionally, the model eliminates the need for a restrictive assumption such as the "mandatory full charge". The replication of each station and battery level based on the number of time intervals provides a direct and intuitive representation of vehicle movement. This approach enables a comprehensive tracking of the overall vehicle status, including their location and battery level. These advantages serve as the primary motivations for adopting the space-time-energy network flow model in formulating the problem addressed in this study.

2.1.2 Integration of vehicle-to-grid with shared electric vehicles

The research on the combination of V2G and SEV is still in its early stages, but its potential and prospects have attracted close attention lately. According to a study conducted by Liao et al. (2021), the integration of V2G services into SEV fleets has been examined for its economic and environmental advantages. The study reveals that, on average, each vehicle participating in V2G services could bring an annual revenue of $2,272 while simultaneously reducing greenhouse gas emissions by 66.5 tons per year. To fully unleash the benefits of V2G and SEV integration, optimizing the charging and discharging processes of SEVs has become a crucial operational aspect. Various methodologies have been introduced to address operational problems stemming from the integration of V2G and SEV.

Freund et al. (2012) is the pioneer study in introducing the integration of V2G with carsharing. They focus on developing an agent-based control architecture for local energy management in a smart distribution feeder. The objective is to maximize energy utilization by considering the involvement of various stakeholders, including the smart grid operator, carsharing operator, and distribution system operator.

Recognizing the significance of price sensitivity in influencing user behavior, Ren et al. (2019) propose a dynamic pricing scheme for a large-scale SEV network integrated with V2G. Using predicted station-level travel demand, they formulate the dynamic pricing problem as an optimization model that considers demand-driven electricity pricing, vehicle relocation, and charging and discharging scheduling.

Considering the slow charging aspect in the context of SEV with V2G operations, Iacobucci et al. (2019) propose an innovative model that optimizes both transport service and charging at distinct time scales using two model-predictive control optimization algorithms. Longer time scales optimize vehicle charging to minimize electricity costs, while shorter time scales focus on vehicle routing and rebalancing to reduce passenger waiting times.

Taking into account uncertainties of carsharing demand and electricity price, Zhang et al. (2021) develop a two-stage stochastic integer programming model to optimize various aspects, including service coverage, EV battery capacity choice, fleet deployment and parking space capacity. To capture the uncertainties, they consider a range of demand-price scenarios. They construct a spatial-temporal-SOC transportation network for each scenario, enabling them to dynamically model daily operations involving EV rentals, relocation, idling, charging, and electricity selling to the grid.

Similarly, considering the stochasticity of electricity price, Li et al. (2022) propose a Markov decision model to determine the optimal time allocations for a fleet of SEVs in both the transportation system and the power market to maximize the total profit. An efficient dynamic programming algorithm is adopted to solve the proposed model. Notably, in their work, the EV fleet is regarded as a whole, without considering detailed routes of individual EVs. This simplification could potentially encounter challenges in practical implementation.

Under the assumption that relocation is performed at non-operation hours, Caggiani et al. (2021) propose optimization models to determine the EV distribution among stations at the beginning of each day to maximize profits from V2G operations while meeting the carsharing demand as much as possible.

Similarly, Prencipe et al. (2022) adopt the static relocation approach and propose a mixed integer linear programming model to determine the EV distribution at the beginning of a day and the optimal daily charging/discharging schedules to maximize the revenue from carsharing users and V2G profits. While static relocation simplifies the problem and reduces computational complexity, its lack of responsiveness to dynamic fluctuations in demand can potentially result in inefficient vehicle utilization.

*2.2 Battery swapping*

Battery swapping, offered by third-party companies or electric vehicle manufacturers enables users to replace their depleted batteries with fully charged ones. This technology has been recognized as an appropriate substitute for conventional fuel stations (Revankar and Kalkhambkar, 2021). Despite its potential, there is a relatively limited amount of literature on battery swapping. Wu (2022a) gives a comprehensive review of the current research on the planning and operational challenges in battery swapping systems. The author highlights that the operational scenarios in battery swapping are intricate, necessitating further advancements in optimization methods. Currently, a commonly used mode is the charging-swapping mode, where users can exchange their batteries for fully charged ones at designated battery swapping stations, and subsequently, their depleted batteries can be recharged at these stations (Shaker et al., 2023). Another mode involves mobile battery swapping, where a battery swapping van delivers fully charged batteries to battery swap demand points and transports depleted ones to battery swapping stations (Yang et al., 2021b). However, our focus is solely on the charging-swapping mode.

Identifying suitable locations and the appropriate number of battery swapping stations is crucial for improving user accessibility and minimizing both construction and operational expenses for operators. The data mining method has been adopted to address this location problem, as demonstrated by Zeng et al. (2019), which proposes a data-driven framework for solving the location selection problem. In their approach, massive GPS data of taxies and an electricity consumption rate model are used to estimate the demand for battery swapping. Alongside data-driven methodologies, various optimization techniques have also been widely employed in this context. For example, when confronted with uncertain demand information, Mak et al. (2013) tackle the battery swapping station deployment issue by formulating two distributionally robust optimization models. The first aims to minimize the robust estimate of anticipated costs, while the second seeks to maximize the worst-case likelihood of attaining a predefined return-on-investment target.

Once the location is obtained, ensuring sustainable operation requires careful consideration of the number of available swappable batteries. Addressing this, Ni et al. (2021) propose a two-step solution scheme for jointly solving the inventory and real-time vehicle-to-station problem. In their approach, the initial inventory is determined using sample average approximation in the first step. Subsequently, based on the optimal inventory solution, a randomized online algorithm in the second step optimizes the real-time vehicle-to-station routing problem. It is important to note that at the strategic-level decision-making, this approach solely addresses inventory planning, with battery swapping station locations regarded as fixed input parameters. Given the high construction cost of battery swapping stations, there is a limited number of stations. Hence, in order to achieve cost-effective operations, it would be more advantageous to integrate both location and inventory considerations.

If the configuration of a battery swapping system is deployed in a region, when an EV requests a swap, the operator assigns a specific station for the driver to swap his/her battery, commonly referred to as the vehicle-to-station assignment problem in the literature. You et al. (2018) study this problem, proposing an optimization model that assigns the most suitable station for swapping a depleted battery to each EV based on its current location and SOC. The assignment aims to minimize a weighted sum of EVs' travel distance and electricity generation cost. Notably, their approach assumes a known quantity of fully charged batteries at stations and does not track charging schedules, which leads to each fully charged battery being replacement for depleted batteries at most once.

With the involvement of microgrids, battery swapping stations can realize B2G interactions (Revankar and Kalkhambkar, 2021; Cui et al., 2023). Through their B2G capability, BSSs may take part in energy and reserve markets and increase their profit. Jordehi et al. (2020) propose a nonlinear programming model to determine the optimal power flow at the batter swapping stations taking into account the B2G operation. Results show the economic benefits with the involvement of B2G. Considering the battery condition, swapping/charging demand, and fluctuating electricity prices over time, the operator faces the task of determining an optimal charging/discharging schedule for the batteries. Confronted with the random customer request problem, Mahoor et al. (2019) develop a robust optimization approach to determine the worst-case solution for the optimal battery involving B2G applications. This schedule includes the charging and discharging power of each battery at each time slot, considering the revenue from power sales to the grid, costs of purchasing power, and battery degradation costs. In their model, the state of each battery at each time is tracked with various states. As the number of batteries increases, the complexity of the problem also grows. However, they do not develop specific algorithms to address this challenge. Consequently, there exist limitations in the practical application of large-scale problems due to the considerable time consumption required.

*2.3 Research gaps*

As previously reviewed, several models and methods have been developed to optimize the recharging process of EVs, but they mostly focus on the plug-in charging. Unfortunately, this recharging method often results in long charging times for EVs, even with fast charging.

In contrast, battery swapping offers a promising alternative that enables EVs recharged fully within a quite short time. It has the potential to complement the plug-in charging approach. It is worth noting that while battery swapping holds potential for EV recharging, there has been no research conducted on its application within carsharing systems, even though carsharing seems to be a suitable context for such an approach.

Recognizing this research gap, we are motivated to study the integration of battery swapping and plug-in charging within carsharing systems. From the previous review, we can observe that the primary concern of models and methods discussed above concentrating on the strategic planning and operational management of battery swapping is ensuring that individual vehicles can access the stations and efficiently exchange their batteries based on the vehicle arrival information. However, the applicability of these models and methods to the context of SEVs is not guaranteed. The integration of battery swapping stations into SEV operations introduces greater complexity. The demand for SEVs and the movement of these vehicles directly influence the need for battery swapping, which subsequently impacts the strategic and operational decisions associated with battery swapping. Therefore, any model or approach aiming to address battery swapping within SEV systems needs to account for this mutual relationship. Surprisingly, to the best of the authors' knowledge, scant research has been conducted in this area, which serves as one of the motivations for this study.

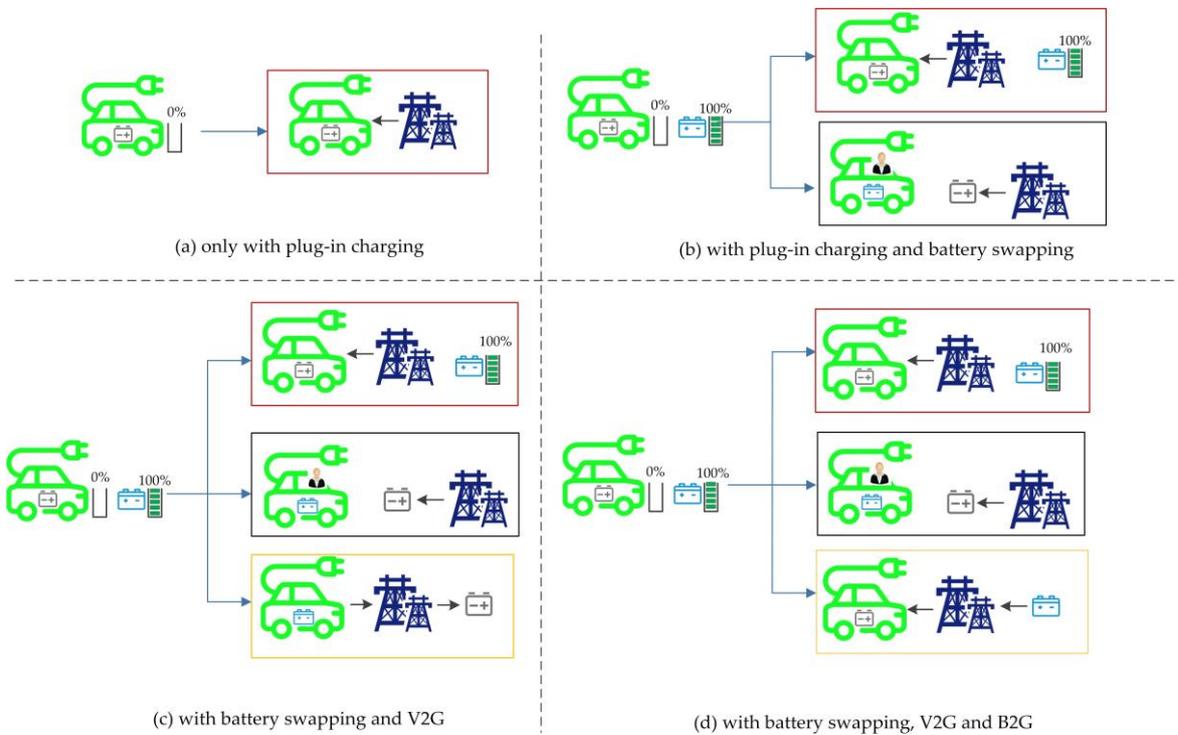

**Fig. 1** Optional strategies under different scenarios.

The integration of various technologies associated with the power energy transmission holds the promise of enhancing satisfaction rate, boosting vehicle and battery utilization, and unlocking additional benefits. To illustrate it, let us examine a scenario depicted in Fig.1. Imagine a user in need of a vehicle. In Fig.1(a), a vehicle with SOC of 0% is present at the station. With only plug-in charging available, this user cannot be served, so the vehicle has to be charged at the station waiting for the next users, as indicated by the red box. In the remaining three scenarios, battery swapping is an option, and a fully charged swappable battery is available alongside the vehicle at the station. Fig. 1(b) illustrates the scenario where both plug-in charging and battery swapping are utilized. Potential strategies involve the vehicle serving the user after a battery swap or opting not to serve the user and instead charging the batteries. In Fig. 1(c), both battery swapping and V2G technology are employed. In addition to the strategies shown in Fig. 1(b), there is the possibility of the vehicle supplying power to the grid following battery swapping due to V2G capabilities. Fig. 1(d) takes it a step further, encompassing battery swapping, V2G, and B2G. In this scenario, there is the potential for the fully charged battery to directly supply power to the grid.

While the potential benefits of integration are evident, existing research, as indicated in Table 1, has typically tackled only one or two of these aspects as optimization problems. However, with each consideration added, the complexity of the problem escalates. More variables come into play when making decisions, and there is a balance to be struck between the costs and revenues associated with each potential strategy. The complexity is also a key reason why a comprehensive optimization framework for charging/discharging scheduling, incorporating the amalgamation of battery swapping, V2G, and B2G, has not been proposed thus far.

To bridge the aforementioned gaps, this study focuses on the optimization of SEV systems the context of V2G and B2G applications, while also considering the integration of plug-in charging and battery swapping.

**Table 1** Summary of previous studies related to the operation of SEVs.

| Articles | Recharging | Decisions | Modeling and solution method | V2G/B2G |
|---|---|---|---|---|
| Xu et al. (2020) | Plug-in charging and mandatory full charge | Charging station deployment | Mixed-integer nonlinear programming model, outer approximation approach | Neither |
| Ait-Ouahmed et al. (2018) | Plug-in charging and mandatory full charge | Relocations | Mixed integer linear programming, tabu search algorithm | Neither |
| Boyaci et al. (2015) | Plug-in charging and mandatory fixed charge | Relocations | Multi-objective stochastic programming | Neither |
| Xu et al. (2018) | Plug-in charging and mandatory full charge | Fleet size, pricing, relocation | Mixed-integer convex programming, outer-approximation algorithm | Neither |
| Zhao et al. (2018) | Plug-in charging and partial charging | Relocation, allocation, charging scheduling | Time-space network, Lagrangian relaxation-based algorithm | Neither |
| Bekli et al. (2021) | Plug-in charging and partial charged | Charger locations, relocation, and assignment | Space-time-energy network, customized heuristic algorithms | Neither |
| Gambella et al. (2018) | Plug-in charging and partial charging | Relocation, assignment, charging | Time-space network, customized algorithm | Neither |
| Huang et al., 2018 | Plug-in charging and partial charging | Fleet size, relocation, charging | Time-space network; customized algorithm | Neither |
| Chen and Liu (2022) | Plug-in charging and partial charging | Location, deployment, fleet size, relocation, charging | Space-time-energy network, two-stage stochastic programming, | Neither |
| Zhang et al. (2019) | Plug-in charging and partial charging | Assignment, relocation, charging | Space-time-energy network, brand and cut algorithm | Neither |
| Gonçalves Duarte Santos et al. (2023) | Plug-in charging and partial charging | Fleet size, charging facilities, relocation, charging | Space-time-energy network | Neither |
| Zhang et al., 2021 | Plug-in charging and partial charging | Location, allocation, relocation, charging/discharging scheduling | Space-time-energy network, Benders decomposition | Only V2G |
| Iacobucci et al. (2019) | Plug-in charging and partial charging | Relocation, charging/discharging scheduling | Model-predictive control optimization algorithm | Only V2G |
| Li et al. (2022) | Plug-in charging and partial charging | Allocation, charging/discharging scheduling | Markov decision model, dynamic programming algorithm | Only V2G |
| Yang et al.(2021b) | Battery swapping | Location, routing | Data-driven approach | Neither |
| This paper | Plug-in charging, | Location, relocation, | Space-time-energy network, CG- | V2G and |

| | partial charging and battery swapping | allocation, charging/discharging scheduling | based heuristic algorithm | B2G |
|---|---|---|---|---|

## 3. Model formulation and small-scale cases

The proposed model pertains to the planning and operation of a reservation-based one-way station-based electric carsharing system. This system allows users to conveniently pick up and return SEVs at designated stations. Moreover, these stations are equipped with charging facilities, and the system has already adopted V2G technology. As referred, the problem we address includes strategic and operational facets. On the strategic front, we determine whether to implement the battery swapping and B2G technology through the upgrade of existing stations. This necessitates assessing which stations warrant an upgrade and determining the optimal quantity and distribution of swappable batteries. At the operational level, our objective is to optimize the charging and recharging of batteries and EVs, as well as EV relocations and battery swapping operations.

Before presenting the mathematical model, the following subsection provides an elaborate description of the model's formulation and definitions that will be used throughout the remainder of this paper.

*3.1 Definitions and formulation*

**Time**. The operational horizon is denoted by $T$, and it is divided into a set of time intervals with uniform durations, denoted as $\{1, 2, \ldots, t, \ldots, |T|\}$, where 1 denotes the beginning of the operational horizon, and $|T|$ denotes the ending time.

**Stations.** Vehicles are picked up and returned at designated stations, denoted as $I = \{1, 2, \ldots, i \ldots, |I|\}$. Each station is equipped with a specific number of parking spaces. Based on the available recharging infrastructure, we consider three types of stations:

- *Type 1*: Stations with only parking spaces, lacking recharging infrastructure.
- *Type 2*: Stations with parking spaces equipped with recharging infrastructure only supporting for V2G technology.
- *Type 3*: Stations with parking spaces featuring battery swapping, as well as V2G and B2G infrastructure. Fig. 2 illustrates the infrastructure layout of a Type 3 station.

In the context of this paper, the construction of Type 3 stations is still pending with Type

3 stations being upgraded from Type 2 stations. Operators strategically determine which Type 2 stations should be upgraded to Type 3. For stations equipped with recharging infrastructure, it is assumed that each parking space is outfitted with a charging outlet. Throughout this paper, we assume that carsharing operators install their own recharging and battery swapping infrastructure.

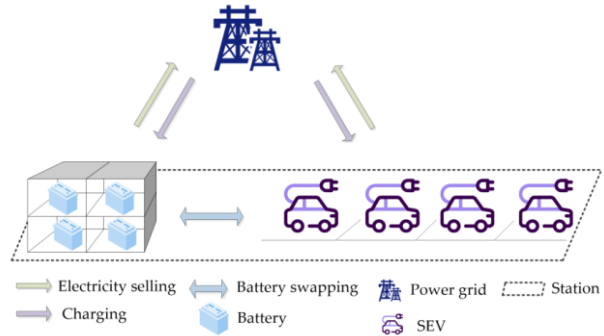

**Fig.2** Illustration of a station.

**Vehicles**. A fixed homogeneous fleet of SEVs provides carsharing service. The operator determines the initial vehicle allocation before the daily operation. Within this system, SEVs can execute a range of operations, outlined as follows:

- *Rental*: Users rent vehicles for one-way trips. Upon receiving a user request, the carsharing operator needs to determine whether the request can be accepted.

- *Relocation*: Due to the demand imbalance, vehicles could accumulate or become scarce at specific stations. Relocation is employed to justify the vehicle distribution. In this study, we consider operator-based relocations, assuming the availability of sufficient staff for performing the relocations. Static overnight relocations during non-operational hours are not considered. The operator should determine the relocation tasks, including when and where vehicle relocations should take place.

- *Charging*: Vehicles can be charged while parked at stations equipped with plug-in charging infrastructure. Partial charging is permissible, meaning vehicles need not be fully charged before their subsequent use. The operator must determine the optimal timing, location, and duration for vehicle charging.

- *Energy Selling*: Beyond serving carsharing users, vehicles can contribute to the grid network by supplying electricity. The operator needs to decide when, where, and for what duration vehicles should sell electricity to the grid using V2G technology.

- *Battery swapping*: Battery swapping involves replacing a vehicle's depleted battery with a fully charged one at dedicated stations. The operator should decide the appropriate timing, location, and number of vehicles for battery swapping operations.

- *Idle*: Vehicles are parked at stations, not in use or actively participating in any of the tasks mentioned above. These vehicles awaiting assignment.

**Charging and discharging**. A nonlinear charging function is adopted, where the battery level increases linearly with charging time at a rate of $\beta_1$ when the battery level is between 0% and 80%. Once the battery level reaches 80%, the charging speed decreases to $\beta_2$. The battery discharging rate is assumed to be linear with the travel time at a constant discharging speed (Chen and Liu, 2022).

**Electricity price**. The electricity price for charging and selling is the same at any given time, but it fluctuates throughout a day.

**Battery**. In this study, batteries owned by operators are uniform and can be classified into two states: *vehicle-installed* and *stocked* (Wang et al., 2019). Stocked batteries are physically kept at Type 3 stations and come with four states: idle, charging, discharging, and swapping. This means that a battery is either installed in a vehicle or stored at a Type 3 station. When the vehicle-installed batteries are nearing depletion, EVs can move to Type 3 stations to swap their empty batteries with the fully-charged stocked ones. When a vehicle-installed battery is swapped out to be stored at a station, simultaneously, another stocked battery is swapped in to become vehicle-installed. At any given time, the total number of either vehicle-installed or stocked batteries remains constant. Specifically, the number of vehicle-installed batteries aligns with the SEV fleet size. Additionally, during the entire operation, stocked batteries in each state transition between states, and the number of stocked batteries stored at each Type 3 station remains constant; however, this exact quantity remains to be determined. At the beginning of a day, it is assumed that all batteries are fully charged. Our study addresses the following sates concerning stocked batteries:

- *Charging*: Typically, at battery swapping stations, each battery is positioned in a charging box before being swapped, thus a removed battery from the vehicle can immediately find its place in a charging box (You et al., 2018).

- *Energy Selling*: With B2G technology, stocked batteries have the capability to directly feed electricity into the grid without relying on electric vehicles. The operator should decide when stocked batteries should sell electricity, taking into account factors such as EV usage demand and electricity prices.

- *Swapping*: Fully charged stocked batteries can be used to replace depleted vehicle-installed batteries.
- *Idle*: In this state, stocked batteries are not engaged in charging, discharging, or swapping activities. They are temporarily unassigned to specific tasks, simply stored at the Type 3 station, awaiting future operations.

Vehicle-installed batteries are installed within vehicles, and their state is synchronized with that of the vehicles they are installed in. In this context, we only detail the operations of stocked batteries.

**Demand**. User requests consist of trip information, including origin, destination, and trip starting time. Demand represents a collection of requests that share common trip characteristics. In this study, demand of a typical day is determined and known.

*3.2 Optimization model based on a space-time-energy network*

A flow-based integer programming model in built on a three-dimensional space-time-energy network, which extends the conventional space-time network flow model by incorporating battery level information. Each node in the network represents a specific combination of location, time, and battery level. Battery levels are discretized into a set of percentage-based levels, denoted by $E$. A node is represented by a triple $v = (i,t,e)$ for $i \in I$, $t \in T$, and $e \in E$. For the sake of model formulation, a virtual source node is introduced, represented as node $(0,0,|E|)$. Here, 0 serves as the dummy time interval and the dummy depot. Arcs between nodes represent the state of vehicles and stocked batteries. Depending on the operations of vehicles and stocked batteries, these arcs can be categorized into the following sets:

**Rental arcs of vehicles**: $(i,t,e,j,k,r) \in A^{rent}$ for $i,j \in I$, $t,k \in T$, and $e,r \in E$, connecting from node $(i,t,e)$ to node $(j,k,r)$ with $k = t + \tau_{ik} \leq |T|$ and $r = e - \beta^D \tau_{ik} \geq 0$, where $\tau_{ik}$ is the shortest travelling time from station $i$ to station $j$. The rental arc indicates the number of vehicles with battery level $e$ which are used by clients to travel from station $i$ to $j$ starting at time instant $t$. For operators, the revenue generated from each unit flow on rental arc $(i,t,e,j,k,r)$ is computed by $p^{rental} \tau_{ij}$, where $p^{rental}$ is the service price of SEV per time interval.

**Relocation arcs of vehicles:** $(i,t,e,j,k,r) \in A^{relo}$ for $i,j \in I$, $t,k \in T$, and $e,r \in E$, with $k = t + \tau_{ik} \leq |T|$ and $r = e - \beta^D \tau_{ik} \geq 0$. The relocation arc indicates the number of vehicles with battery level $e$ which are relocated by staff to travel from station $i$ to $j$ beginning at time instant $t$. For operators, the unit cost on the relocation arc $(i,t,e,j,k,r)$ is $c_{ij}^{relo}$.

**Idle arcs of vehicles:** $(i,t,e,i,t+1,e) \in A^{idle}$ for $i \in I$, $1 \leq t \leq |T|-1$, and $e \in E$, flows on these arcs represent the number of EVs remaining idle at station $i$ between time interval $t$ and time interval $t+1$ with SOC equal to $e$. In this case, these vehicles do not change either in space or energy but only ahead in time.

**Charging arcs of vehicles:** $(i,t,e,i,t+1,r) \in A^{charg}$ for $i \in I_2, 1 \leq t \leq |T|-1$, and $0 \leq e \leq |E|-1$, with $r = e + \beta_e^C \leq |E|$. Flows on these arcs represent the number of EVs being charged from $e$ to $r$ at station $i$ between time interval $t$ and time interval $t+1$. For operators, the unit cost on charging arc $(i,t,e,i,t+1,r)$ is $\lambda(r-e)p_t^{elec}$, where $\lambda$ is the maximum capacity of a battery, expressed in kWh, and $p_t^{ele}$ is the electricity price at time interval $t$.

**Energy selling arcs of vehicles:** $(i,t,e,i,t+1,r) \in A^{sell}$ for $i \in I_2, 1 \leq t \leq |T|-1$, $0 < e \leq |E|$ with $r = e - \beta^D \geq 0$. Flows on these arcs represent the number of EVs selling electricity back to the grid at station $i$ between time interval $t$ and time interval $t+1$ and the SOC decreasing from $e$ to $r$. For operators, the revenue generated from each unit flow on the energy selling arc $(i,t,e,i,t+1,r)$ is $\lambda(e-r)p_t^{elec}$.

**Battery swapping arcs of vehicles:** $(i,t,0,i,t+1,|E|) \in A^{swap}$ for $i \in I_2$, $1 \leq t \leq |T|-1$, it is assumed that it takes one time step to complete battery swapping. Flows on these arcs represent the number of EVs with SOC of $e$ swapping batteries at station $i$ between time interval $t$ and time interval $t+1$. For operators, the unit cost on a battery swapping arc is $c^{swap}$.

**Idle arcs of stocked batteries:** $(i,t,e,i,t+1,e) \in B^{idle}$ for $i \in I_2$, $1 \leq t \leq |T|-1$, and $e \in E$. Flows on these arcs represent the number of stocked batteries with SOC of $e$ idling at station $i$ between time interval $t$ and time interval $t+1$.

**Charging arcs of stocked batteries:** $(i,t,e,i,t+1,r) \in B^{charg}$ for $i \in I_2$, $1 \leq t \leq |T|-1$ and $0 \leq e \leq |E|-1$, with $r = e + \beta_e^C \leq |E|$. Flows on these arcs represent the number of stocked batteries being charged from $e$ to $r$ at station $i$ between time interval $t$ and time interval $t+1$. For operators, the cost per unit flow on charging arc $(i,t,e,i,t+1,r)$ is $\lambda(r-e)p_t^{elec}$.

**Energy selling arcs of stocked batteries:** $(i,t,e,i,t+1,r) \in B^{sell}$ for $i \in I_2$, $t \in [1,|T|-1]$, $0 < e \leq |E|$ with $r = e - \beta^D \geq 0$. Flows on these arcs represent the number of stocked batteries selling electricity back to the grid at station $i$ between time interval $t$ and time interval $t+1$ with SOC decreasing from $e$ to $r$. For operators, the revenue generated from each unit flow on the energy selling arc $(i,t,e,i,t+1,r)$ is $\lambda(e-r)p_t^{elec}$.

**Battery swapping arcs of stocked batteries:** $(i,t,|E|,i,t+1,0) \in B^{swap}$, for $i \in I_2$, $t \in [1,|T|-1]$. Flows on these arcs represent the number of stocked batteries being swapped with empty inside-vehicle batterie at station $i$ between time interval $t$ and time interval $t+1$.

**Dummy source arcs:** $(0,0,|E|,i,1,|E|) \in A^{source}$ for $i \in I$, $(0,0,|E|,i,1,|E|) \in B^{source}$ for $i \in I_2$. These arcs indicate initial allocation of EVs and stocked batteries to station $i$ at the start.

The notations used in our formulation are summarized in Table 2.

**Table 2** Notations for the model formulation.

| | Sets and index |
|---|---|
| $I$ | Set of stations, indexed by $i$ and $j$ |
| $I_1$ | Set of stations solely used for parking without recharging infrastructure |
| $I_2$ | Set of stations with recharging infrastructure only supporting V2G |
| $I_3$ | Set of battery swapping stations with recharging infrastructure supporting V2G and B2G |
| $T$ | Set of time intervals, indexed by $t$ and $k$ |
| $E$ | Set of battery SOC, indexed by $e$ and $r$ |
| $A^{rent}$ | Set of vehicle rental arcs in the space-time-energy network |
| $A^{relo}$ | Set of vehicle relocation arcs in the space-time-energy network |
| $A^{idle}$ | Set of vehicle idle arcs in the space-time-energy network |
| $A^{charg}$ | Set of vehicle charging arcs in the space-time-energy network |
| $A^{sell}$ | Set of vehicle selling energy arcs in the space-time-energy network |
| $A^{swap}$ | Set of swapping batteries arcs of vehicles in the space-time-energy network |
| $A^{source}$ | Set of vehicle source arcs in the space-time-energy network |
| $A$ | The union of sets related to vehicle activities, that is, $A = A^{rent} \cup A^{relo} \cup A^{idle} \cup A^{charg} \cup A^{sell} \cup A^{swap} \cup A^{source}$ |
| $A^{park}$ | A subset of $A$, $A^{park} = A^{idle} \cup A^{charg} \cup A^{sell} \cup A^{swap} \cup A^{source}$ |
| $B^{idle}$ | Set of stocked battery idle arcs in the space-time-energy network |
| $B^{charg}$ | Set of stocked battery charging arcs in the space-time-energy |
| | network |
| $B^{sell}$ | Set of stocked battery selling energy arcs in the space-time-energy network |
| $B^{swap}$ | Set of stocked battery swapping arcs in the space-time-energy network |
| $B^{source}$ | Set of stocked battery source arcs in the space-time-energy network |
| $B$ | The union of sets related to stocked battery activities, that is, $B = B^{idle} \cup B^{charg} \cup B^{sell} \cup B^{swap} \cup B^{source}$ |
| $\zeta_{ite}^{+}$, $\zeta_{ite}^{-}$ | Set of arcs for which $(i,t,e)$ is the origin and destination node of EVs, respectively |
| $\theta_{ite}^{+}$, $\theta_{ite}^{-}$ | Set of arcs for which $(i,t,e)$ is the origin and destination node of stocked batteries, respectively |
| $A_{it}^{start}$ | A subset of $A$, arcs originating from station $i$ starting from time instant $t$, $A_{it}^{start} = \{(i,t,e,j,k,r) \in A, e \in E, j \in I, k \in T, r \in E\}$ |
| $A_{it}^{swap}$ | A subset of $A^{swap}$, swapping arcs from station $i$ to station $j$ starting from time instant $t$, $A_{it}^{swap} = \{(i,t,0,i,t+1,|E|) \in A^{swap}\}$ |
| $B_{it}^{swap}$ | A subset of $B^{relo}$, swapping arcs from station $i$ to station $j$ starting from time instant $t$, $B_{it}^{swap} = \{(i,t,|E|,i,t+1,0) \in B^{swap}\}$ |
| $A_{ijt}^{relo}$ | A subset of $A^{relo}$, relocation arcs from station $i$ to station $j$ starting from time instant $t$, $A_{ijt}^{relo} = \{(i,t,e,j,k,r) \in A^{relo}, e \in E\}$ |
| $A_{ijt}^{rent}$ | A subset of $A^{rent}$, rental arcs from station $i$ to station $j$ starting from time instant $t$, $A_{ijt}^{rent} = \{(i,t,e,j,k,r) \in A^{rent}, e \in E\}$ |
| Parameters ||
| $p^{rental}$ | Revenue from a SEV per time interval when rented by a client |
| $p_t^{ele}$ | Unit electricity prices of charging or selling at time interval $t$ |
| $c_{ij}^{relo}$ | Cost of relocating a SEV from $i$ to $j$ |
| $c^{swap}$ | Cost of swapping a battery for a SEV |
| $c^{battery}$ | Depreciation cost of a battery per day |
| $c^{station}$ | Daily cost of upgrading a Type 2 station to Type 3 |
| $p_a$, $p_b$ | Revenue from arc $a$, $b$ |
| $c_a$, $c_b$ | Cost of arc $a$, $b$ |
| $\lambda$ | Maximum capacity of a battery, expressed in kWh |
| $\tau_{ij}$ | Travel time between stations $i$ and $j$ |
| $n_i$ | Number of parking spaces of station $i$ |
| $m_i$ | Capacity for storing stocked batteries at station $i$ |
| $\beta_e^C$ | Charging speed at battery level $e$, in percentage |
| $\beta^D$ | Discharging speed at battery level $e$, in percentage |
| Decision variables ||

| | | |
|---|---|---|
| $z$ | | Number of stocked batteries |
| $x_a$ | | Flow of vehicles on arc $a \in A$ |
| $y_b$ | | Flow of stocked batteries on arc $b \in B$ |
| $s_i$ | | Binary variable, 1 if station $i \in I_2$ is chosen for the upgrade to Type 3, and 0 otherwise |

The integrated optimization framework is formulated as a space-time-energy network flow model, denoted as [P]. The objective and constraints of this model are as follows.

$$\max \sum_{a \in A^{rent} \cup A^{sell}} p_a x_a + \sum_{b \in B^{sell}} p_b y_b - \sum_{a \in A^{relo} \cup A^{charg} \cup A^{swap}} c_a x_a - \sum_{b \in B^{charg} \cup B^{swap}} c_b y_b - c^{battery} z - c^{station} \sum_{i \in I_2} s_i \quad (1)$$

$$\text{s.t.} \sum_{a \in A^{source}} x_a = F \quad (2)$$

$$\sum_{a \in \zeta_{ite}^+} x_a - \sum_{a \in \zeta_{ite}^-} x_a = 0 \quad \forall i \in I, t \in \{1,2,\ldots,|T|-1\}, e \in E \quad (3)$$

$$\sum_{a \in A_{ijt}^{rent}} x_a \leq d_{ijt} \quad \forall i,j \in I, t \in T \quad (4)$$

$$\sum_{a \in \zeta_{ite}^- \cap A^{rent}} x_a \leq \sum_{a \in \zeta_{ite}^+ \cap A^{idle}} x_a \quad \forall i \in I, t \in T, e \in E \quad (5)$$

$$\sum_{a \in A^{start} \cap A^{park}} x_a \leq n_i \quad \forall i \in I, t \in T \quad (6)$$

$$\sum_{b \in B^{source}} y_b = z \quad (7)$$

$$\sum_{b \in \theta_{ite}^+} y_b - \sum_{b \in \theta_{ite}^-} y_b = 0 \quad \forall i \in I_2, t \in \{1,2,\ldots,|T|-1\}, e \in E \quad (8)$$

$$y_b \leq m_i s_i \quad \forall i \in I_2, b = (0,0,|E|,i,1,|E|) \in B^{Source} \quad (9)$$

$$\sum_{a \in A_{it}^{swap}} x_a = \sum_{b \in B_{it}^{swap}} y_b \quad \forall i \in I_2, t \in T \quad (10)$$

$$x_a, y_b, z \in \mathbb{Z} \quad \forall a \in A, b \in B \quad (11)$$

$$s_i \in \{0,1\} \quad \forall i \in I_2 \quad (12)$$

Objective function (1) maximizes the overall profit, which is a function encompassing both revenues and costs. Revenues arise from fulfilled user requests and the sale of electricity from both SEVs and batteries. On the other hand, costs consist of the relocation costs, electricity buying costs, battery swapping costs, stocked battery depreciation costs and daily depreciation station upgrade cost.

Constraint (2) ensures that the total initial number of EVs assigned to stations is equal to the fleet size. This is grounded on the assumption that all vehicles originate from the dummy depot. Notably, fleet size $F$ is considered a parameter rather than a decision variable in this study.

Constraints (3) represent the flow conservation constraints at space-time-energy nodes. Constraints (4) require that the number of EVs assigned to the trip requests with origin $i$, destination $j$ and departure time $t$ should not exceed the demand for the OD pair. Constraints (5) ensure that there are enough number of idle vehicles to start trips at each station at each time interval. Constraints (6) ensure that the number of vehicles parked at station $i$ at the beginning of a given time interval cannot exceed the capacity of the station. It should be noted that for station $i \in I_1$, the parked vehicles include idle vehicles; for station $i \in I_2$, the parked vehicles include idle vehicles, vehicles being charged and vehicles selling electricity; and for station $i \in I_3$, the parked vehicles include idle vehicles, vehicles being charged, vehicles selling electricity and vehicles that are swapping their batteries.

Constraint (7) calculates the total number of stocked batteries needed, which also indicate their allocation at stations. Constraints (8) ensure the flow conservation for stocked batteries. The purpose of constraints (9) is two-fold. They both prevent batteries stocked at stations which are not upgraded and restrict the stocked battery allocation not exceed the station capacity for batteries.

Constraints (10) indicate that the battery swapping process involves two simultaneous actions. For EVs, an empty inside-vehicle battery is removed, while a stocked battery comes. For batteries, an empty inside-vehicle battery transforms into a stocked one, while a fully charged battery becomes an inside-vehicle one. Constraints (11)-(12) define the domains of decision variables $x_a$, $y_b$, $z$ and $s_i$, specifying that $x_a$, $y_b$ and $z$ are integers, while $s_i$ is a binary variable.

*3.3 Small-scale illustration*

We design a small-scale network to demonstrate the output of the proposed model, the test network has the following characteristics: there are 4 stations, 8 time intervals, 3 SEVs, and 6 requests. Station 1 is categorized as Type 1, while the remaining belongs to Type 2. SOC level is discretized into the following levels: {0%, 10%, 20%, 30%, 40%, 50%, 60%, 70%, 80%, 90%, 100%}. When the SOC is between 0% and 80%, the charging rate per time interval is 40%, while beyond 80%, the charging speed diminishes to 10% per time interval. Meanwhile, the discharging rate remains constant at 10% per time interval. Additionally, the electricity price is higher during time interval 3 to 4 compared to other intervals. The travel time between any two stations is displayed in Fig. 3, and it should be noted that the travel time is symmetric. Trip demand information including the origin, destination, departure time, and quantity is shown in Table 3. The solution of the test instance is depicted in Fig. 3. In this figure, x-axis represents time, y-axis represents SOC levels, z-axis represents the station index, and each space-time-energy node is denoted by a grey circle. Colored arrows represent arcs, while the numbers on the arcs represent the flow on the arcs. Fig. 4(a) illustrates the SEV flows, while Fig. 4(b) shows the battery flows.

As depicted in Fig. 4(a), at the beginning, two EVs are assigned to Station 1, while the remaining EV is allocated to Station 4. All requests are satisfied. One relocation is performed from Station 4 to Station 1 beginning at time instant 4. Throughout this relocation, the vehicle's SOC drops from 60% to 30%, as indicated by the yellow arrow. Two vehicles sell electricity from time interval 3 to 4, and one vehicle sells electricity from time interval 6 to time interval 8, as shown by black arrows. To serve another trip, a battery swapping operation is conducted at Station 2 during time interval 5 and time interval 6, as depicted by the blue arrow.

Fig. 4(b) displays the battery flow, and it is observed that Station 2 is upgraded, and one stocked battery is required. Taking advantage of the higher electricity price during time intervals 3 to 4, this fully charged battery participates in electricity sales. Following this, it is recharged to its maximum capacity. From time interval 5 to time interval 6, this battery is placed into a vehicle, thereby substituting the battery with a SOC of 0%. The battery with 0% SOC remains at Station 2, remaining idle until the end of time interval 8.

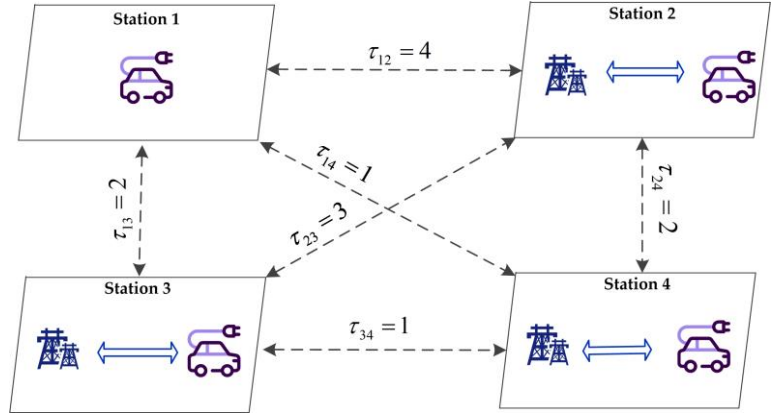

**Fig. 3.** Travel time between stations (in number of time intervals)

**Table 3** Trip demand information.

| Origin | Destination | Departure time | Quantity |
|--------|-------------|----------------|----------|
| 1 | 2 | 1 | 1 |
| 2 | 3 | 5 | 1 |
| 1 | 4 | 1 | 1 |
| 3 | 4 | 4 | 1 |
| 4 | 3 | 1 | 1 |
| 1 | 4 | 5 | 1 |

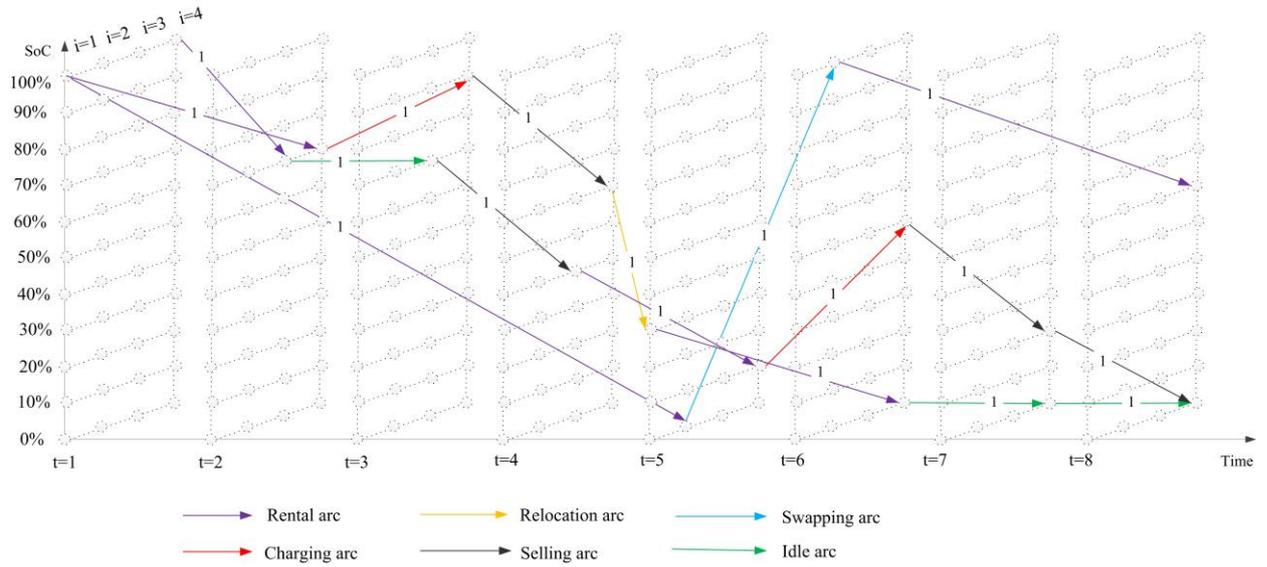

(a)

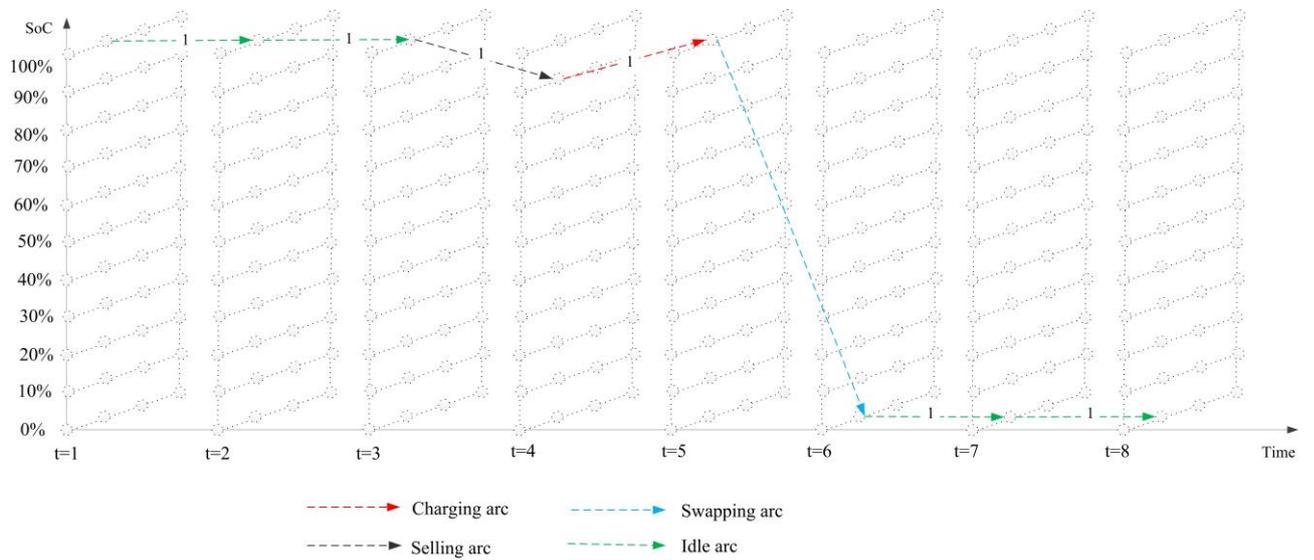

(b)

**Fig. 4.** Flow results of the illustration case.

## 4. Solution approach

*4.1 Computation complexity and model preprocessing*

With an extra dimension for tracking battery levels, the state space of the proposed network flow model can expand significantly as the problem scale grows larger. Furthermore, the considered problem is integrated, as it involves the interconnection of vehicle flows and battery flows. The mathematical model [P] has a set of decision variables ($x_a, y_b, s_i, z$), and Table 4 summarizes the number of decision variables for each set, enabling an analysis of the model complexity. Here, $|\bullet|$ denotes the cardinality of set $\{\bullet\}$.

**Table 4** Number of decision variables per set.

| Decision variables | Number of variables |
| --- | --- |
| $x_a$ | $2|I|^2|T||E|+|I||T||E|+2|I_2||T||E|+|I||T|+|I|$ |
| $y_b$ | $2|I_2||T||E|+2|I_2||T|+|I_2|$ |
| $s_i$ | $|I_2|$ |
| $z$ | 1 |

In the network model, the number of decision variables depends on the number of arcs. To reduce the size of problem, we conduct preprocess to only consider the necessary arcs and nodes. Among the various sets of arc types, rental arcs and relocation arcs have the largest number of arcs, with a number of $|I|^2|T||E|$, which makes it time-consuming to obtain a solution. Hence, it becomes imperative to preprocess them to generate the minimum number of arcs while preserving essential connectivity.

To minimize the number of rental arcs, we generate rental arcs $A_{ijt}^{rent}$ only if a user request exists from station $i$ to $j$ at time $t$. If there is no demand from station $i$ to $j$ at time $t$, then $A_{ijt}^{rent}$ is set to be empty. Regarding the relocation arcs, we adopt the approach introduced in Chen et al. (2022) of identifying hub stations to reduce the number of relocation arcs. For the relocation trip between station $i$ and $j$, we seek a hub $h$ where the travel time adheres to the condition $\tau_{ij} = \tau_{ih} + \tau_{hj}$. If such a hub $h$ exists, the relocation arcs between station $i$ and $j$ can be eliminated. This is because the relocation trip from station $i$ to $j$ can be completed by using two separate relocation trips: one from station $i$ to $h$ and another from station $h$ to $j$.

*4.2 Column generation-based heuristic algorithm*

4.2.1 Reformulation

As introduced, this study considers four activities of stocked batteries: idle, charging, selling energy, and battery swapping. Besides, we adopt the model setting detailed in the work of You et al., 2018. According to their settings, at Type 3 stations, each stocked battery is placed in a charging box. When a depleted inside-vehicle battery is swapped with a fully charged stocked battery, the depleted battery will be deposited in the same box that previously held the fully charged one. This setting indicates that if a box in the beginning of the day is set to have a battery in it, then there will always be a battery in this box, even as the battery itself and its activities may change. Conversely, if a box is not initially assigned any batteries, it will remain empty throughout the day. Therefore, we can make the following observation:

**Observation 1**: The number of stocked batteries is equal to the number of used boxes. By tracking the state of used boxes during each time interval, the flow of stocked batteries is obtained.

Fig. 5 illustrates this observation by considering a station equipped with three distinctively colored boxes across four time intervals. Only one stocked battery is assigned to this station and placed in the first box. Throughout the first time interval, the battery in the first box remains idle. During the second time interval, it sells energy. The third time interval sees it charging, and in the fourth time interval, a battery swap is executed. Clearly, at this station, only one box is used, leaving the remaining two boxes empty. It can be observed that the number of used boxes aligns with the initially allocated number of stocked batteries. This observation holds for any given station. Furthermore, within the entire system, the number of stocked batteries is equal to the number of used boxes.

The state of a box is defined to be consistent with the activity of the battery positioned inside it. As an illustration, in Fig. 5, the state of the first box is categorized as idle during the first time interval, because the battery inside it is idle. To illustrate the relationship between the box states and the flow of stocked batteries, we introduce state arcs for each box. For the sake of simplification, we will solely consider the idle state arcs as an example in this context. Initially, we construct the idle state arcs for boxes. Following the process of arc construction for batteries, we define idle state arcs $(i',t,e,i',t+1,e) \in B_{box}^{idle}$ for box $i'$, which means that there is an idle stocked battery with an SOC $e$ placed in box $i'$ during time interval between $t$ and $t+1$. Subsequently, we define a binary variable $v_{b'}$ ($b'=(i',t,e,i',t+1,e)\in B_{box}^{idle}$), which takes the value 1 if there is an idle stocked battery on arc $b'$, and 0 otherwise. Then, the stocked battery flow $y_b$ on arc $b=(i,t,e,i,t+1,e)\in B^{idle}$ can be expressed as $y_b = \sum_{b'} v_{b'}$ ($b'\in\{(i',t,e,i',t+1,e)\in B_{box}^{idle}\,|\,i'\in I_i^{box}\}$), where $I_i^{box}$ is the set of boxes at the station $i$. This equation signifies that the flow of idle stocked batteries at station $i$ is the cumulative sum of boxes belonging to station $i$ and containing idle stocked batteries. By extending this logic, various other types of stocked battery flows can also be represented using the state arcs associated with boxes. This aligns with the statement made in Observation 1, where by tracking the state of boxes during each time interval, the flow of stocked batteries is obtained.

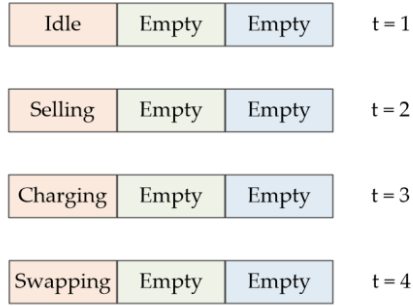

Fig. 5. States of boxes.

Based on Observation 1, Constraints (7)-(10) can be replaced with the following constraints:

$$\sum_{b'\in B_{box}^{source}} v_{b'} = z \tag{13}$$

$$\sum_{b'\in \vartheta_{i'te}^+} v_{b'} - \sum_{b'\in \vartheta_{i'te}^-} v_{b'} = 0 \quad \forall i'\in I^{box}, t\in\{1,2,\ldots,|T|-1\}, e\in E \tag{14}$$

$$v_{b'} \leq s_i \quad \forall i\in I_2, b'\in\{(0,0,|E|,i',1,|E|)\in B_{box}^{source}\,|\,i'\in I_i^{box}\} \tag{15}$$

$$\sum_{a\in A_{it}^{swap}} x_a = \sum_{b'\in B_{it}^{swap\_box}} v_{b'} \quad \forall i\in I_2, t\in T \tag{16}$$

$$v_{b'}\in\{0,1\} \quad \forall b'\in B^{box} \tag{17}$$

where $I^{box}$ is the set of all boxes, $B^{box}$ is the set of all state arcs of boxes, and $B_{box}^{source}$ is the dummy source arcs of boxes. $\vartheta_{i'te}^+$ and $\vartheta_{i'te}^-$ are a set of arcs for which $(i',t,e)$ is the origin and destination node, respectively. $B_{box}^{swap}$ is the set of swapping arcs of boxes, and $B_{it}^{swap\_box}$ is a subset of $B_{box}^{swap}$, denoted by $B_{it}^{swap\_box} = \{(i',t,|E|,i',t+1,0) \in B_{box}^{swap} | i' \in I_i^{box}\}$. Other notations correspond to those introduced in Section 3.2. By simply substituting the station index in the stocked battery arcs with the box index, all box state arcs can be generated. Hence, the detailed enumeration of these box state arcs is not provided here. Constraints (13)-(14), and (16) are equivalent of Constraints (7)-(8) and (10), respectively, but written in terms of the variable $v_{b'}$ instead of $y_b$. Constraints (15) ensure that only the station $i$ is selected to upgrade to Type 3 station, then stocked batteries can be placed within the boxes in station $i$.

Compared to the formulation in Section 3.2, the above reformulation leads to heightened computational complexity due to the increased number of decision variables. The total number of decision variables has now increased by $(2|I_2||T||E| + 2|I_2||T| + |I_2|) \cdot \left(\sum_{i \in I_2} m_i - 1\right)$. However, this reformulation has an advantage, it is expected to yield solutions by using CG algorithm. In the next subsections, we explain how to apply the CG algorithm to our proposed model.

4.2.2 Master problem

CG is a specialized technique that aims to address large-scale optimization problems. This approach has demonstrated in successful applications in diverse areas such as routing problems (Wang et al., 2023; He et al., 2023), scheduling problems (Zhou et al., 2022; Wu et al., 2022b), and more. The fundamental concept of CG is to iteratively solve a master linear programming model and a pricing problem. The master linear programming model includes a subset of variables of the original problem, and it passes dual variables to the pricing problem to find promising columns with negative reduced costs for minimization problems (or positive reduced costs for maximization problems). Subsequently, the decision variables that correspond to these promising columns are incorporated into the master linear programming model. This iterative approach continues until no promising columns can be generated.

For any given box, a sequence of its state arcs within a day can be represented as a **chain**. This chain consists of a series of state arcs in an ascending order in terms of time. Take the first box in Fig. 5 as an example, suppose that the power of the battery inside the box remains at 100% from time interval $t=1$ to $t=2$, then changes from 100% to 60% during time interval $t=2$ to $t=3$, subsequently returns from 60% to 100% during time interval $t=3$ to $t=4$, and finally drops from 100% to 0 during time interval $t=4$ to $t=5$. The chain representation for this box across these time intervals is as follows: $(1,100\%,2,100\%)_{idle} \Rightarrow (2,100\%,3,60\%)_{selling} \Rightarrow (3,60\%,4,100\%)_{charging} \Rightarrow (4,100\%,5,0)_{swapping}$. Box chains provide a concise reflection of the stocked battery operations. Thus, determining the optimal schedule for stocked batteries equates to finding the most favorable box chains. Considering the pool of feasible box chains available, the formulation presented in Sections 4.2.1 and 3.2 can be stated as the following master problem ([MP]):

$$\max \sum_{a \in A^{rent} \cup A^{sell}} p_a x_a - \sum_{a \in A^{relo} \cup A^{charg}} c_a x_a - c^{battery} \sum_{i \in I_2, u \in U} g_{iu} - c^{station} \sum_{i \in I_2} s_i + \sum_{i \in I_2, u \in U} c_u^{chain} g_{iu} \quad (18)$$

s.t. Constraints (2), (4)-(5)

$$\sum_{a \in \zeta_{ite}^+ \setminus A^{swap}} x_a + \sum_{u \in U, a \in \zeta_{ite}^+ \cap A^{swap}} \delta_{iu}^a g_{iu} - \sum_{a \in \zeta_{ite}^- \setminus A^{swap}} x_a - \sum_{u \in U, a \in \zeta_{ite}^- \cap A^{swap}} \delta_{iu}^a g_{iu} = 0 \quad \forall i \in I, t \in \{1, 2, \ldots, |T|-1\}, e \in E$$

(19)

$$\sum_{a \in A_{it}^{start} \cap A^{park} \setminus A^{swap}} x_a + \sum_{u \in U, a=(i,t,0,i,t+1,|E|) \in A^{swap}} \delta_{iu}^a g_{iu} \leq n_i \quad \forall i \in I, t \in T \quad (20)$$

$$\sum_{u \in U} g_{iu} \leq m_i s_i \quad \forall i \in I_2 \quad (21)$$

$$x_a, y_b, g_{iu} \in \mathbb{Z} \quad \forall a \in A, b \in B, i \in I_2, u \in U \quad (22)$$

$$s_i \in \{0,1\} \quad \forall i \in I_2 \quad (23)$$

where $U$ represents the set of feasible chains, and $\delta_{iu}^a$ is an indicator parameter signifying the presence of a swapping arc within chain $u$ at station $i$. When chain $u$ at station $i$ contains a swapping arc $(t,|E|,t+1,0)$, it indicates the occurrence of a vehicle swapping its battery during time interval $t$ to $t+1$. Hence, we let $\delta_u^a = 1$ $(a=(i,t,0,i,t+1,|E|))$ if chain $u$ contains the arc $(t,|E|,t+1,0)$, and 0 otherwise. $c_u^{chain}$ denotes the profit from chain $u$. The parameters are derived from the solution of the pricing subproblem, which will be presented in the following subsection. $g_{iu}$ is the decision variable, presenting the number of boxes selecting chain $u$, with these boxes belonging to station $i$. It should be noted that each used box matches a single chain, but one station could be linked with more than one chain when it accommodates multiple boxes. Hence, the decision variable $g_{iu}$ is not a binary variable but an integer variable. Furthermore, it is worth mentioning that different boxes can potentially share the same chain, implying that the state of these boxes remains equal throughout the day.

The objective function continues to have the same meaning as in [P]. Each selected chain corresponds to a used box, and according to Observation 1, the number of used boxes equals the number of stocked batteries, so $\sum_{i \in I_2, u \in U} g_{iu}$ represents not only the number of selected chains but also the number of stocked batteries. Hence, the term $c^{battery} \sum_{i \in I_2, u \in U} g_{iu}$ is the depreciation associated with purchasing stocked batteries. The term $\sum_{i \in I_2, u \in U} c_u^{chain} g_{iu}$ is equal to $\sum_{b \in B^{sell}} p_b y_b - \sum_{b \in B^{charg} \cup B^{swap}} c_b y_b$. Constraints (19)-(21) are equivalent to Constraints (3), (6)-(10), respectively, but expressed in terms of the variable $g_{iu}$ instead of $y_b$.

While [MP] is completely equivalent to [P], as the problem size increases, the number of feasible chains might experience exponential growth, creating a significant challenge for solving. Moreover, it becomes impossible to enumerate all feasible chains. To address this concern, CG deals with the restricted master problem ([RMP]), which involves a subset of feasible chains, enabling a more manageable computation. Within the iteration process of CG, new promising columns are generated using information from the dual variables. To achieve this, the relaxation of [RMP] ([RMLP]) is tackled. In this relaxed version, Constraints (22)-(23) are relaxed to continuous variables ($x_a, y_b, g_{iu} \geq 0$, $0 \leq s_i \leq 1$). This relaxation enhances the exploration of potential solutions and contributes to the computation efficiency.

4.2.3 Pricing problem

To formulate the pricing problem, an expanded notation is introduced, as shown in Table 5.

**Table 5** Additional variables and parameters for the pricing problem.

| | Sets and index |
|---|---|
| $H^{sell}$ | Set of selling energy state arcs of a box, representing the state whereby a box |
| | contains a battery that is selling energy, $H^{sell} = \{(t,e,t+1,e-\beta^D) \mid t \in T, e \in E\}$ |
| $H^{charg}$ | Set of charging state arcs of a box, representing the state whereby a box contains a box that is being charged $H^{charg} = \{(t,e,t+1,e+\beta_e^C) \mid t \in T, e \in E\}$ |
| $H^{idle}$ | Set of idle state arcs of a box, representing the state whereby a box contains a battery that is idle $H^{idle} = \{(t,e,t+1,e) \mid t \in T, e \in E\}$ |
| $H^{swap}$ | Set of swapping battery state arcs of a box, representing the state whereby a box contains a battery that is replacing a depleted vehicle-stalled battery $H^{swap} = \{(t,|E|,t+1,0) \mid t \in T\}$ |
| $H$ | Set of all state arcs, $H = H^{sell} \cup H^{charg} \cup H^{idle} \cup H^{swap}$ |
| $\Phi_{te}^+, \Phi_{te}^-$ | Set of state arcs for which $(t,e)$ is the origin and destination nodes, respectively |
| Parameters | |
| $\pi_{ite}^1, \pi_{it}^2, \pi_i^3$ | Dual values |
| $p_h$ | Revenue from state arc $h$, with a definition consistent with $p_b$ in Section 3.2 |
| $c_h$ | Cost of arc $h$, with a definition consistent with $c_b$ in Section 3.2 |
| Variables | |
| $f_h$ | Binary variable, 1 if state arc $h$ is selected to compose a chain, and 0 otherwise |

After solving [RMLP], the dual values for Constraints (19)-(21) can be obtained. Let $\pi_{ite}^1$ represent the dual values of Constraints (19), $\pi_{it}^2$ represent the dual values of Constraints (20), and $\pi_i^3$ represent the dual values of Constraints (21). For each station $i$, the reduced cost of column $u$ in [RMLP] can be formulated as follows:

$$\overline{c}_{iu} = c_u^{chain} - \sum_{t \in T, e \in E} \pi_{ite}^1 ( \sum_{a \in \zeta_{ite}^+ \cap A^{swap}} \delta_u^a - \sum_{a \in \zeta_{ite}^- \cap A^{swap}} \delta_u^a) + \sum_{t \in T} \pi_{it}^2 \sum_{a=(i,t,0,i,t+1,|E|) \in A^{swap}} \delta_u^a + \pi_i^3 \quad (24)$$

where

$$c_u^{chain} = \sum_{h \in H^{sell}} p_h f_h - \sum_{h \in H^{charg} \cup H^{swap}} c_h f_h \quad (25)$$

By performing the necessary substitutions and collecting terms, this reduced cost can be expressed as follows:

$$\overline{c}_{iu} = \sum_{h \in H^{sell}} p_h f_h - \sum_{h \in H^{charg} \cup H^{swap}} c_h f_h - \sum_{t \in T} \sum_{h=(i,t,|E|,i,t+1,0) \in H^{swap}} (\pi_{it|E|}^1 - \pi_{it0}^1) f_h + \sum_{t \in T} \pi_{it}^2 \sum_{h \in H^{swap}} f_h + \pi_i^3 \quad (26)$$

In each iteration of the CG process, $|I_2|$ pricing subproblems are solved, with each subproblem associated with a specific station $i \in I_2$. A feasible solution to a pricing subproblem corresponds not only to a feasible chain of a box, but also to an element within the set $U$ of [MP]. The pricing problem is to find a promising column with a positive reduced cost, which subsequently is added to [RMLP]. The formulation of the subproblem for each station can be presented as an integer linear program, as detailed below. It should be noted that in each subproblem, the index $i \in I_2$ is removed from all decision variables.

[SP$_i$]:

$$\max \sum_{h \in H^{sell}} p_h f_h - \sum_{h \in H^{charg} \cup H^{swap}} c_h f_h - \sum_{t \in T} \sum_{h=(t,|E|,t+1,0) \in H^{swap}} \left( \pi^1_{it|E|} - \pi^1_{it0} \right) f_h + \sum_{t \in T} \pi^2_{it} \sum_{h \in H^{swap}} f_h + \pi^3_i \quad (27)$$

s.t.
$$\sum_{h \in \Phi_{1|E|}} f_h = 1 \quad (28)$$

$$\sum_{h \in \Phi^+_{te}} f_h = \sum_{h \in \Phi^-_{te}} f_h \quad \forall t \in T, e \in E \quad (29)$$

$$\sum_{h \in H^{swap}} f_h \geq 1 \quad (30)$$

$$f_h \in \{0,1\} \quad \forall h \in H \quad (31)$$

Constrain (28) restricts a box to have only one type of state during the initial time interval. Constraints (29) ensure the conservation of states of a box at node $(t,e)$. Constraints (29) and (30) together ensure that a box only has one state arc within any time interval. Constraint (30) guarantees that at least one battery swap happens within a box. Chains without swapping arcs are straightforward to generate—simply organize charging and selling arcs based on the principle of selling electricity when electricity prices are high. Within the process of our algorithm, chains without swapping arcs are added into the columns of [RMLP] during the initialization phase. Hence, Constraint (30) is incorporated into the model to prevent the generation of chains without swapping arcs. Constraints (31) impose binary conditions on the variable $f_h$. The non-zero solutions of $f_h$, arranged in ascending order of time, form a chain of boxes at station $i$.

Since each subproblem is solved separately for each station, the complexity of the problem is significantly reduced. The termination criterion for the CG iteration is met when the largest computed reduced cost, denoted as $\bar{c}^* = \max_i \bar{c}_{iu}$, is negative. In general, after terminating the CG algorithm, if the current solution to [RMLP] is non-integer, the branch and bound method should be incorporated into the CG technique, known as the branch and price algorithm. Our experimentation with this algorithm has faced challenges, primarily in the selection of effective branching rules due to the high dimensionality of the decision variables. The absence of effective branching rules leads to a large execution time. Consequently, we directly solve the [RMP] with the current column set (Du et al., 2014). After solving the [RMP], we obtain the final near-optimal solution of the original problem. The proposed algorithm sometimes can find the optimal solution. When solving the [RMP], if the chains selected in the optimal solution of the original problem are also part of the current chain set, the proposed algorithm can discover optimal solutions to the original problem. The framework of the proposed CG-based heuristic algorithm is shown in Fig. 6.

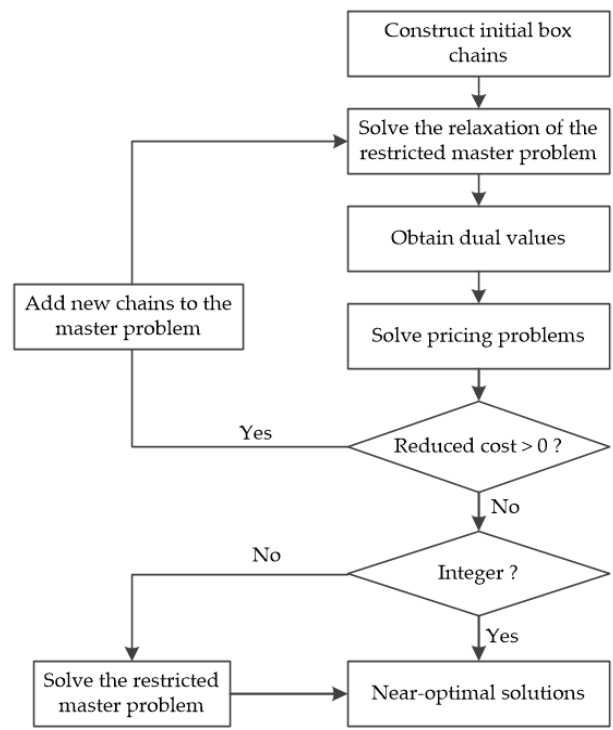

**Fig. 6**. Framework of the CG-based heuristic algorithm.

4.2.4 Accelerating techniques

**Column elimination.** When the number of columns in [RMLP] becomes too large, a subset of the non-basic columns can be removed. Specifically, columns that have not entered the basis of the master problem in the last $\eta$ iterations (where $\eta$ is a predefined parameter) are considered for removal. To ensure proper convergence of the CG process, this elimination procedure should not be performed frequently, and a certain number of non-basic columns should remain in [RMLP] after column removal.

**Diversity**. Column diversity aids in comprehensively exploring the solution space, thereby increasing the chances of finding the optimal solution. In this study, we use the number of swaps in a chain to describe the diversity. When solving the pricing subproblem, Constraint (30) is substituted with the subsequent constraint:

$$\sum_{h \in H^{swap}} f_h \geq \kappa \tag{32}$$

where $\kappa$ is a parameter. The value of this parameter ranges from 1 to $\lfloor |T|/\alpha \rfloor + 1$, where $\alpha$ is the time required for a battery to be continuously charged from empty to full, and $\lfloor |T|/\alpha \rfloor + 1$ represents the maximum number of battery swaps that occur at a box. We assign various values to the parameter $\kappa$, with each value corresponding to a distinct subproblem model. These models can be solved in parallel.

**Multiple columns.** To reduce the number of iterations, at each iteration, we add more than one column with a positive reduced cost (if any). It is not necessary to find the chain with the maximal reduced cost, as long as the reduced cost of the chain is greater than 0, the objective function of [RMLP] can be improved. In the computational experiments, we use the solution pool feature provided by CPLEX solver, a widely utilized commercial software for solving mathematical optimization problems, to collect multiple solutions, whether they are optimal or not, generated throughout the optimization process when solving [SP$_i$]. However, it should also be noted that, if too many columns are added to each iteration, the size of [RMLP] will become increasingly larger, subsequently extending the solution time accordingly.

**Relaxation.** To reduce computation time, we opt to solve a mixed integer programming model ([MIP]) instead of solving the [RMP] to obtain the final near-optimal solution. The only difference between this [MIP] and the [RMP] lies in the domain of decision variables. In the [RMP], all variables are integers. However, in the [MIP], we only specify the variables $x_a (a \in A^{rent} \cup A^{source})$, $g_{iu}$ and $s_i$ as integers, while the others are relaxed, as shown below.

[MIP]:

$$\max \sum_{a \in A^{rent} \cup A^{sell}} p_a x_a - \sum_{a \in A^{relo} \cup A^{charg}} c_a x_a - c^{battery} \sum_{i \in I_2, u \in U} g_{iu} - c^{station} \sum_{i \in I_2} s_i + \sum_{i \in I_2, u \in U} c_u^{chain} g_{iu} \tag{33}$$

s.t. Constraints (2), (4)-(5), (19)-(21),

$$x_a \geq 0 \quad \forall a \in A \tag{34}$$

$$x_a, y_b, g_{iu} \in \mathbb{Z} \quad \forall a \in A^{rent} \cup A^{source}, b \in B, i \in I_2, u \in U \tag{35}$$

$$s_i \in \{0,1\} \quad \forall i \in I_2 \tag{36}$$

It is worth noting that the optimal solutions of the [RMP] and the [MIP] are the same, as Proposition 1 outlines, with its proof provided in Appendix A.

**Proposition 1.** There exists an optimal solution to the [MIP] where all decision variables are integer.

**Proof.** See Appendix A.

## 5. Numerical experiments

In Section 5.1, we test a set of numerical instances to validate the proposed algorithm. Furthermore, in Section 5.2, we perform a series of sensitivity analyses based on a real-world case in the city of Lanzhou, Gansu, China, to identify key parameters and assess their impact on the system's performance. All experiments are carried out using a personal laptop equipped with an Intel Core i5 processor running at 2.4GHz and 16 GB of RAM.5.2

*5.1 Computation performance of the solution method*

In this section, we assess the performance of the proposed CG-based heuristic algorithm through 13 randomly generated examples. The effectiveness of the algorithm is compared against CPLEX. Our computations are executed under a predefined maximum computation time of 10800 seconds. If CPLEX surpasses this time limit, the computation is terminated, and results from CPLEX are reported as "-". The summarized results are presented in Table 6.

In Table 6, Column 1 denotes the index of problems, Column 2 denotes the problem scale represented by a triple (# of stations, # of time intervals, # of demand), Column 3 denotes the optimal objective value attained by CPLEX, Column 4 indicates the computation time required by CPLEX to find the optimal solution, Column 5 is the objective value found by the CG-based heuristic algorithm, Column 6 reports the time taken to obtain the solution using the CG-based heuristic algorithm, and Column 7 is the objective value gap between the solutions obtained by CPLEX and the CG-based heuristic algorithm, which is defined as follows:

$$\text{Obj. gap} = \frac{\text{Column 3-Column 5}}{\text{Column 3}} \times 100\% \tag{37}$$

Analyzing the data in Table 6, we observe that both solution approaches require increased computation time as the instance scale grows. Notably, the performance of CPLEX solutions deteriorates more rapidly compared to that of the CG-based heuristic algorithm. Instances 10-13, pose a challenge for CPLEX, as it fails to discover the optimal solution within 10800s. For relatively smaller-scale instances, such as Instances 1, 2, and 3, both methods yield the optimal solution swiftly. In the case of the preceding 9 instances, where CPLEX identifies the optimal solutions, the objective value gap between these two methods ranges from 0% to 3.51% with an average of 1.53%. As the scale increases, the distinct advantage of the CG-based heuristic algorithm becomes increasingly evident. Our proposed algorithm outperforms CPLEX in terms of running time for Instances 5-13. Following these results, the proposed algorithm is able to find good solutions within a reasonable time and its performance notably surpasses that of the commercial solver for slightly larger instances.

**Table 6** Computation performance of the two solution approaches.

| Instances | Problem scale | CPLEX | | CG-based heuristic algorithm | | Obj. gap (%) |
|---|---|---|---|---|---|---|
| | | Obj. | Time(s) | Obj. | Time (s) | |
| 1 | (10, 10, 100) | 850 | 0.59 | 850 | 2.13 | 0 |
| 2 | (10, 15, 200) | 2501 | 1.56 | 2501 | 5.04 | 0 |
| 3 | (10, 20, 300) | 5074 | 13.22 | 5074 | 19.12 | 0 |
| 4 | (20, 15, 400) | 5241 | 22.47 | 5201 | 34.52 | 0.76 |
| 5 | (20, 20, 300) | 6138 | 63.18 | 6021 | 50.16 | 1.90 |
| 6 | (20, 30, 500) | 11662 | 220.02 | 11495 | 91.32 | 1.43 |
| 7 | (25, 30, 600) | 15759 | 761.93 | 15204 | 305.07 | 3.51 |
| 8 | (30, 40, 1000) | 30065 | 2349.66 | 29183 | 864.51 | 2.93 |
| 9 | (50, 30, 1500) | 37214 | 9113.89 | 36002 | 3075.21 | 3.25 |
| 10 | (60, 30, 1500) | - | - | 39552 | 4208.56 | - |
| 11 | (55, 35, 2000) | - | - | 50932 | 5761.17 | - |
| 12 | (50, 40, 2000) | - | - | 52412 | 6384.03 | - |
| 13 | (60, 35, 2000) | - | - | 53671 | 6929.65 | - |

*5.2 Case study and sensitivity analysis*

Exiangxing, founded in 2017 in the city of Lanzhou, China, is station-based car-sharing company that provides one-way and two-way car-sharing services. In the subsequent case studies, we employ the actual operational data from Exiangxing to evaluate the efficacy of the proposed method. The original database covers the time period of August 2018. For the analysis purpose, we select the top 40 stations with the highest demand during this particular month, as illustrated in the layout depicted in Fig. 7. Furthermore, the cumulative demand for these 40 stations throughout August is illustrated in Fig. 8a, while Fig. 8b s displays the demand for a typical day. Additionally, Fig. 8c shows the OD demand between the 40 stations for each hour of the day. Notably, the demand distribution is imbalanced, with peak periods occurring in the morning and evening. We pick 7am to 7pm as the study period of the case study, which covers the peak period and the nonpeak period. It is assumed that the relocation time, driving time, and distance between stations are calculated using the shortest path. SOC level is discretized into the following levels: {0%, 10%, 20%, 30%, 40%, 50%, 60%, 70%, 80%, 90%, 100%}. The parameters used in the case study are based on the Exiangxing operation data, and their settings are detailed in Table 7.

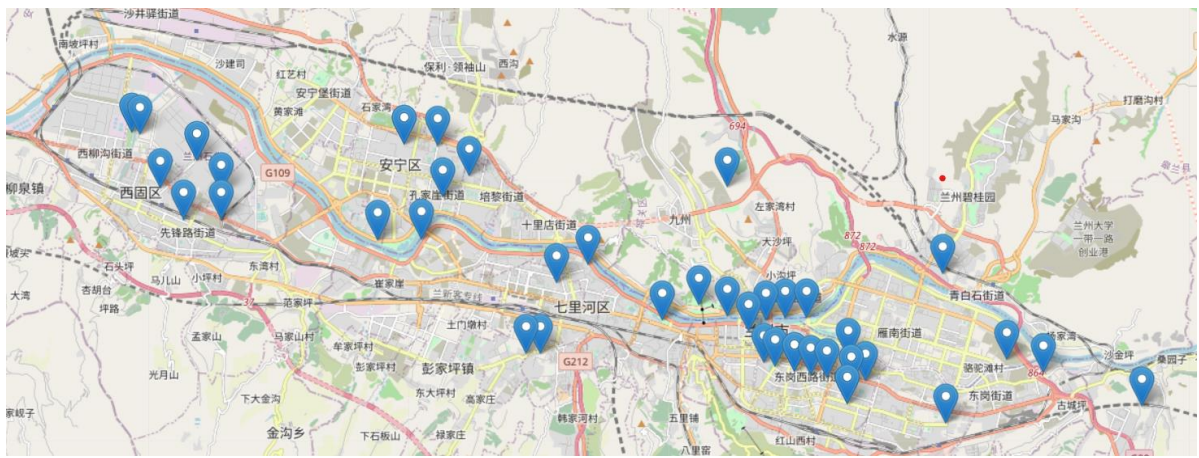

**Fig. 7.** Exiangxing station locations in the city of Lanzhou.

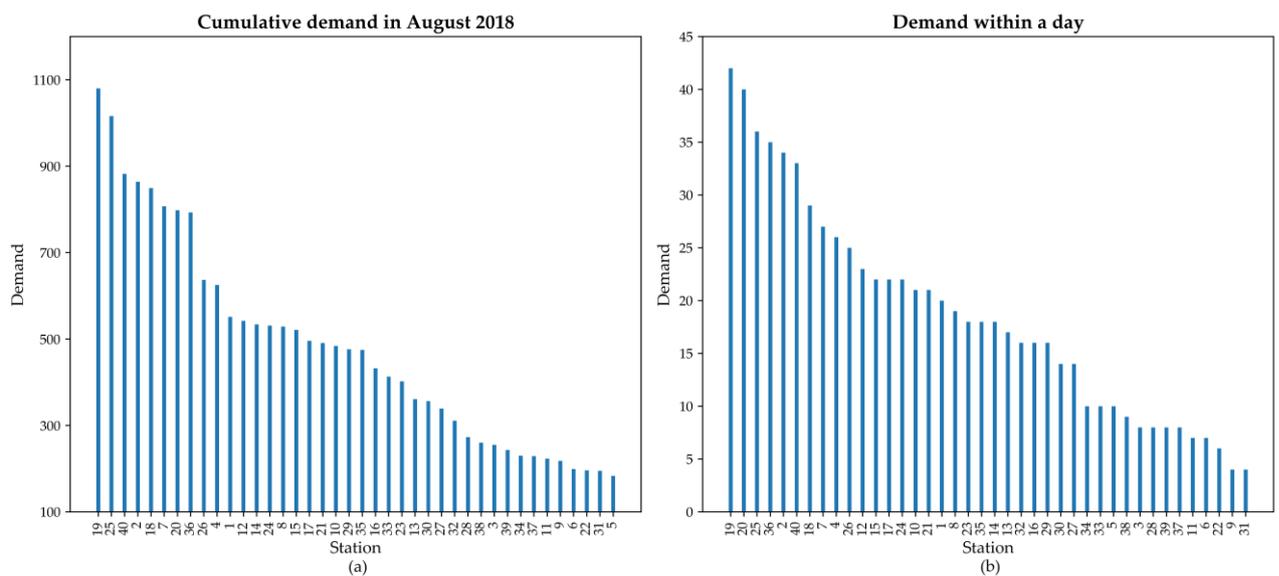

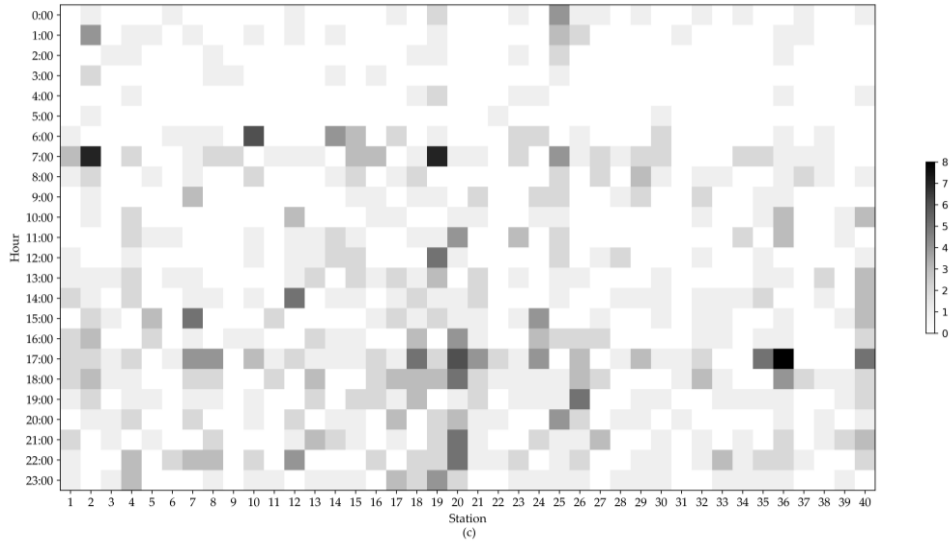

**Fig. 8.** Monthly, daily and hourly demand.

**Tabel 7** Parameter settings in the case study.

| | |
|---|---|
| Number of time intervals | 52 |
| Number of stations | 40 |
| Length of a time interval | 15 min |
| Fleet size | 100 |
| Amount of demand | 1500 |
| Planning horizon | 7:00am – 19:00pm |
| Electricity price in peak periods | 0.759 yuan/kWh |
| Electricity price in off-peak periods | 0.261 yuan/kWh |
| Electricity price in flat-peak periods | 0.51 yuan/kWh |
| Service price | 0.8 yuan/min |
| Relocation cost | 0.3 yuan/min |
| Swapping cost | 5 yuan |
| Number of parking spaces in a station | 5 |
| Capacity of a locker | 5 |
| Depreciation cost of a battery per day | 15 yuan/day |
| Upgrade cost of a station per day | 25 yuan/day |

5.2.1 Sensitivity analysis on the charging speed

In this section, we analyze the impact of charging speed on the system, as charging speed directly affects the operational efficiency of the vehicles. Different charging rates can potentially lead to variations in system performance, making it necessary to conduct in-depth investigations. In line with this requirement, we perform sensitivity analysis across three different charging rates: (i) fast charging, (ii) normal charging, and (iii) slow charging. In the case of fast charging and normal charging, the charging curve follows a non-linear pattern. Up to an SOC of 80%, fast charging operates at a rate of 40%, while normal charging maintains a rate of 20%. On the other hand, the slow charging curve is linear and operates at a constant rate of 10%. These charging profiles are depicted in Fig. 9. Results are reported in Table 8 and Fig.10.

We can observe from Table 8 that as the charging speed increases, both profits and the number of satisfied requests increase. This is attributed to the reduced charging duration, which affords more time to serve passengers. Table 8 shows that the proportion of time allocated to passenger service is higher in scenarios with faster charging.

An interesting observation is that the charging duration proportion of fast charging is 0.34% lower compared to normal charging, but it is 0.43% higher than that of slow charging. When comparing normal charging and fast charging, the time taken for charging of fast charging is shorter, resulting in a lower duration proportion.

However, the speed difference between slow charging and fast charging is significant. Slow charging involves more battery swapping operations, as shown in Table 8, resulting in a higher number of stocked batteries. As a result, for slow charging, a portion of the charging time is transferred to the batteries. The proportion of charging duration for slow charging is 21.89% higher than fast charging. Consequently, the charging duration proportion for vehicles utilizing slow charging is lower than that of fast charging.

Fig. 10 illustrates the quantities of charging and selling energy for vehicles and stocked batteries across different hours. The grey bars in the figure represent electricity prices, which fluctuate throughout the study period. Fig. 10a and Fig. 10b present the total charging and selling energy quantities in the system. Fig. 10c and Fig. 10d show the quantities of charging and selling energy quantities for vehicles, and Fig. 10e and Fig. 10f depict the corresponding quantities for stocked batteries. The values shown in Fig. 10c and Fig. 10e for each hour sum up to the values in Fig. 10a, while the values in Fig. 10d and Fig. 10f correspondingly sum up to the values in Fig. 10b. It is evident that during peak electricity price periods, both vehicles and stocked batteries tend to avoid purchasing energy and instead lean towards selling energy. This indicates that V2G and B2G technologies indeed play a role in energy storage and load balancing. Since the number of batteries in fast charging and normal charging modes is fewer than in the slow charging mode, charging and selling energy quantities are higher for batteries under slow charging, as shown in Fig. 10e and Fig. 10f. An interesting observation is that despite the lower proportion of charging duration and selling energy duration for EVs in the fast charging mode compared to the normal charging mode, the fast charging mode results in higher quantities of energy purchased and energy sold, as depicted in Fig. 10c and Fig. 10d.

Another noteworthy phenomenon is that the electric vehicles do not participate much in the process of selling electricity. The amount of electricity sold by electric vehicles is less than one-tenth of that sold by stocked batteries. The revenue generated by electric vehicles primarily depends on serving carsharing demands. Despite the high electricity prices at 7am and 6pm, which are also peak hours for carsharing, the electric vehicles prioritize serving passengers. Therefore, the amount of electricity sold by the vehicles during these peak hours is not high. This indicates that the emergence of V2G and B2G does not affect the core business of catering to the travel demands of carsharing services.

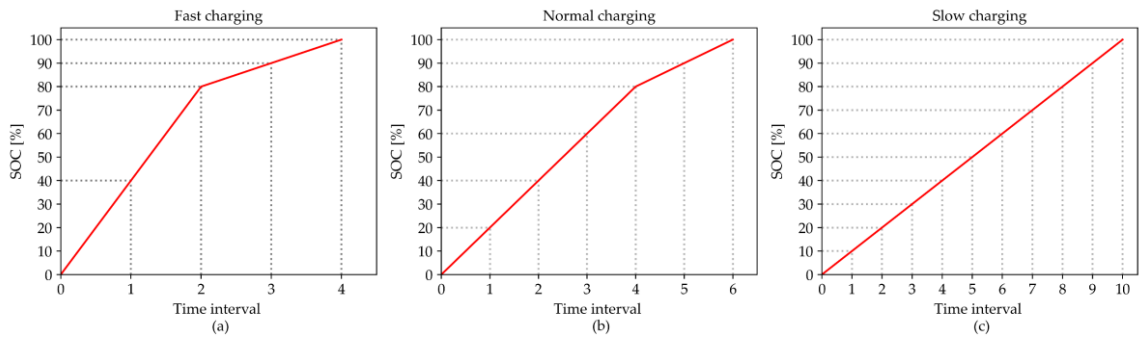

**Fig. 9**. Charging curves.

**Table 8** Sensitivity analysis results using different charging speeds.

|  | Fast charging | Normal charging | Slow charging |
|---|---|---|---|
| Profit | 16954.72 | 16592.71 | 16036.35 |
| Number of satisfied requests | 950 | 936 | 908 |
| Operations of battery swapping | 10 | 51 | 82 |
| Number of stocked batteries | 4 | 17 | 30 |
| Number of Type 3 stations | 2 | 6 | 11 |
| Vehicle time: | | | |
| Moving users (%) | 51.82% | 51.04% | 50.67% |
| Relocation time (%) | 2.71% | 2.48% | 2.07% |
| Charging time (%) | 8.89% | 9.23% | 8.46% |
| Selling energy time (%) | 0.74% | 0.41% | 0.34% |

| Battery time: | | | |
|---|---|---|---|
| Charging time (%) | 21.74% | 32.18% | 43.63% |
| Selling energy time (%) | 45.29% | 37.66% | 31.04% |

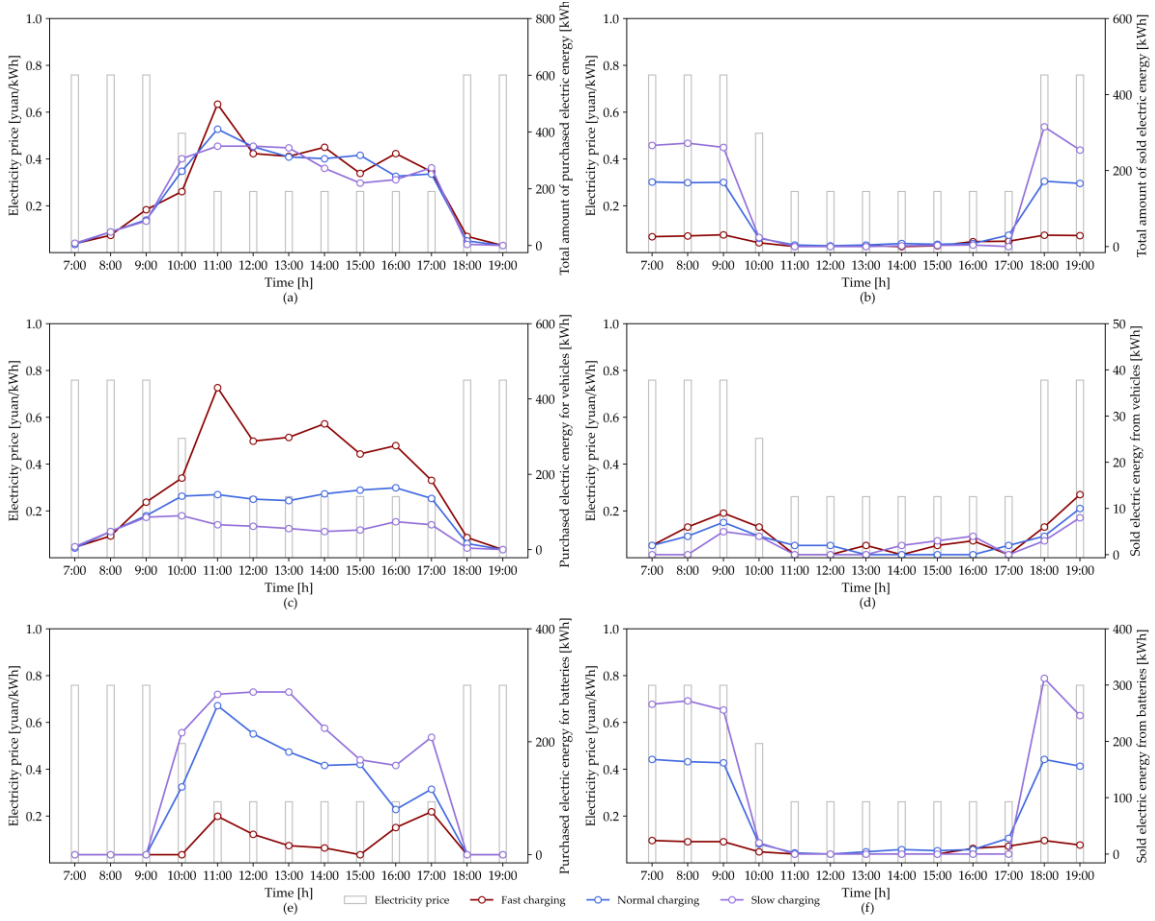

**Fig. 10.** Charging and selling energy quantity each hour.

5.2.2 Sensitivity analysis on the fleet size

In this section, we conduct a sensitivity analysis on the changes in the number of vehicles. In our experimental setup, we consider a range of fleet sizes from 40 to 150 vehicles, with increments of 10. The resulting indicators are presented in Fig. 11. As the fleet size increases, there are corresponding changes in the depreciation costs of the vehicles. This leads to a modification in the objective function of [P], taking the following form:

$$\sum_{a \in A^{rent} \cup A^{sell} \cup A^{sell}} p_a x_a + \sum_{b \in B^{sell}} p_b y_b - \sum_{a \in A^{relo} \cup A^{charg}} c_a x_a - \sum_{b \in B^{charg} \cup B^{swap}} c_b y_b - c^{battery} z - c^{station} \sum_{i \in I_2} s_i - c^{vehicle} F \quad (38)$$

where $c^{vehicle}$ is the depreciation cost of a vehicle, $F$ remains a parameter.

Fig.11a illustrates the variation in profit and average profit per vehicle as the fleet size increases. While profit grows with the increase of the fleet size, the rate of profit increase diminishes, leading to a decrease in average profit per vehicle. Fig. 11b depicts the change in service rate, indicating the increase in served demands as the fleet size grows. However, beyond a fleet size of 100, the rate of increase slows down noticeably. Furthermore, it is observed that the average travel duration for served demands decreases. Particularly, between fleet sizes of 40 and 80 vehicles, the travel duration shows the most significant decrease, from 31.33 to 26.79, a reduction of 14.49%. This suggests that with fewer vehicles that cannot sufficiently meet the demands, prioritizing the demand with longer travel durations. Fig. 11c shows the variations in vehicle relocation time, vehicle charging time, vehicle selling energy time, outside-battery charging time, and outside-battery selling energy time as the fleet size changes. As the fleet size increases, relocation time referring to the time taken on performing relocations, shows a rapid growth trend. Initially, due to the limited vehicle supply relative to demand, relocation is scarcely required, and the system can function primarily based on vehicle borrowing and returning. As the fleet size grows, accumulation of vehicles at certain stations becomes more prevalent, necessitating more relocation tasks.

Fig. 11d demonstrates the variations in the number of battery swaps, the number of stocked battery quantity, and the number of battery swapping stations as the fleet size changes. It is worth noting that as the fleet size increases, these three performance indicators initially rise and then decline. With a smaller fleet size, serves as a means to swiftly recharge vehicles to meet as much as possible demand, compensating for the limited number of vehicles. Between fleet sizes of 40 and 70, the rise in battery quantity aligns with the more gradual growth in the vehicle count. In this range, the demand for battery swapping remains essential due to the insufficient supply of vehicles compared to the demand for car-sharing service, resulting in a higher number of battery swaps and stocked batteries. However, once the fleet size surpasses 70, the need for battery swapping diminishes. This finding further validates that battery swapping is an auxiliary charging method primarily suitable for scenarios with fewer available vehicles. As the fleet size increases, the dependency on battery swapping diminishes. With a sufficiently large fleet, the adoption of battery swapping becomes unnecessary.

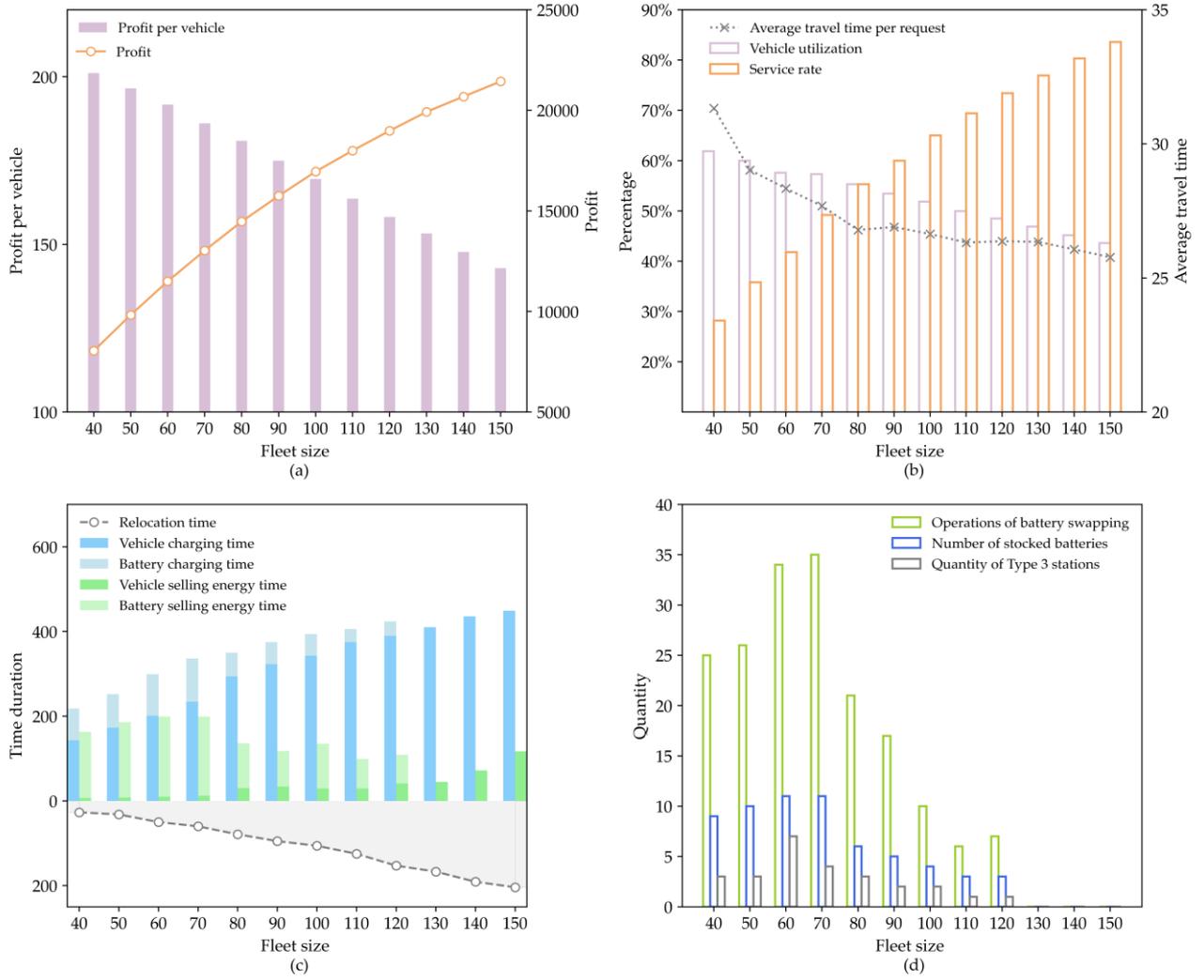

**Fig.11.** Results of sensitivity analysis of fleet size.

## 6. Conclusion

In this paper, we propose an integer programming model based on the space-time-energy network to address various strategic and operational aspects of one-way station-based electric carsharing systems with V2G and B2G integration. The model includes battery swapping operations, station upgrades, charging and selling energy, as well as vehicle relocations. Through a CG-based heuristic algorithm, we efficiently tackle the complex optimization problem presented. To validate our proposed approach, we conduct experiments using both toy and real-world cases. Our findings offer insights into the significance of key parameters within SEV system operations. For example, by comparing fast, normal, and slow charging scenarios, we find that fast charging yields higher profits, and slow charging is more reliant on battery swapping. Moreover, it is observed that while a larger fleet can enhance service rates and profits, it can also lead to reduced vehicle utilization rates. This requires operators to carefully balance operational efficiency and profitability when determining the fleet size. Additionally, in cases with fewer vehicles, there is a greater dependency on battery swapping technology.

Several potential aspects for future research could enhance and extend the scope of this study. Firstly, it would be interesting to incorporate the relocation of stocked batteries between swapping stations into the problem formulation. Secondly, integrating stochastic or uncertain elements inherent to carsharing systems would yield a more realistic and robust solution. Moreover, while this study optimizes charging and discharging schedules, it does not factor in battery degradation. Thus, developing strategies that mitigate battery degradation during charging and discharging is a direction for further investigation.

**Acknowledgements**

The first author wishes to thank the Department of Transport and Planning, Delft University of Technology, Delft, The Netherlands, for hosting her in 2022-2023. Thanks go also to the Chinese Scholarship Council for sponsoring the first author. This work is supported by the National Natural Science Foundation of China [grant numbers 72288101, 71890972/71890970, 72361137003].

**Appendix A. Proof of the Proposition 1**

This proof relies on two theorems related to the totally unimodular matrix:

**Theorem 1**: Each collection of columns of a matrix can be split into two parts so that the sum of the columns in one part minus the sum of the columns in the other part is a vector with entries only 0, 1 and -1, this matrix is totally unimodular (Camion, 1965).

**Theorem 2**: If the coefficient matrix of a linear program is a totally unimodular matrix, and the right-hand side values of the constraints are integers, then the integrality property holds for the linear program (Schrijver, 1986).

For any given combination of variables ($x_{a(a \in A^{rent} \cup A^{source})}$, $g_{iu}$, $s_i$), assume this combination can yield feasible solutions for the [MIP]. Now, the remaining decision variables in the [MIP] are $x_{a(a \in A^{charg} \cup A^{sell} \cup A^{idle} \cup A^{relo})}$. For simplicity, we denote $A_1 = A^{charg} \cup A^{sell} \cup A^{idle} \cup A^{relo}$, denote $A_2 = A^{charg} \cup A^{sell} \cup A^{idle}$, denote $\sum_{a \in A^{rent}} p_a x_a - c^{battery} \sum_{i \in I_2, u \in U} g_{iu} - c^{station} \sum_{i \in I_2} s_i + \sum_{i \in I_2, u \in U} c_u^{chain} g_{iu}$ as $\psi$, and denote $\sum_{a \in \zeta_{ite}^{+} \cap A^{rent}} x_a + \sum_{u \in U, a \in \zeta_{ite}^{+} \cap A^{swap}} \delta_{iu}^{a} g_{iu} - \sum_{a \in \zeta_{ite}^{-} \cap A^{rent}} x_a - \sum_{u \in U, a \in \zeta_{ite}^{-} \cap A^{swap}} \delta_{iu}^{a} g_{iu}$ as $\gamma_{ite}$. We consider the following problem:

$$L_{(x_{a(a \in A^{rent} \cup A^{source})}, g_{iu}, s_i)} = \max \sum_{a \in A^{sell}} p_a x_a - \sum_{a \in A^{relo} \cup A^{charg}} c_a x_a + \psi \tag{39}$$

s.t. $\sum_{a \in \zeta_{ite}^{+} \cap A_1} x_a - \sum_{a \in \zeta_{ite}^{-} \cap A_1} x_a = -\gamma_{ite} \quad \forall i \in I, t \in \{1, 2, \ldots, |T|-1\}, e \in E \tag{40}$

$\sum_{a \in A_{it}^{start} \cap A_2} x_a \leq n_i - \sum_{u \in U, a=(i,t,0,i,t+1,|E|) \in A^{swap}} \delta_{iu}^{a} g_{iu} \quad \forall i \in I, t \in T \tag{41}$

$x_a \geq 0 \quad \forall a \in A_1 \tag{42}$

Given specific values of ($x_{a(a \in A^{rent} \cup A^{source})}$, $g_{iu}$, $s_i$), then $\psi$, $\gamma_{ite}$, and $n_i - \sum_{u \in U, a=(i,t,0,i,t+1,|E|) \in A^{swap}} \delta_{iu}^{a} g_{iu}$ are fixed and must be integers. The coefficient matrix of the above model can be represented as $\begin{bmatrix} \mathbf{C} \\ \mathbf{D} \ \mathbf{0} \end{bmatrix}$, where $\mathbf{C}$ corresponds to the coefficients of variables $x_{a(a \in A_1)}$ in Constraints (40), $[\mathbf{D} \ \mathbf{0}]$ corresponds to the coefficients of variables $x_{a(a \in A_1)}$ in Constraints (41), $\mathbf{D}$ represents the coefficients of $x_{a(a \in A_2)}$ in Constraints (41), and $\mathbf{0}$ represents the coefficients of $x_{a(a \in A^{relo})}$ in Constraints (41). $\mathbf{C}$ can be regarded as the adjacency matrix of the minimum most network flow problem, in which each column contains a single 1 and a single -1, with all other elements being 0. For $\mathbf{D}$, each column contains a single 1, with all other elements being 0. The coefficient matrix in this model can be divided into two parts: $\mathbf{C}$ corresponds to the first part, and $[\mathbf{D} \ \mathbf{0}]$ corresponds to the second part. The sum of each column in one part minus the sum of each column in the other part is a vector with entries only 0, and -1. According to **Theorem 1**, this matrix is totally unimodular. As per **Theorem 2**, the optimal solution for this model is guaranteed to be integers. Each combination ($x_{a(a \in A^{rent} \cup A^{source})}$, $g_{iu}$, $s_i$) corresponds to a $L_{(x_{a(a \in A^{rent} \cup A^{source})}, g_{iu}, s_i)}$. Solving the [MIP] is equivalent to find $(x_{a(a \in A)}, g_{iu}, s_i) \in \arg\max \{L_{(x_{a(a \in A^{rent} \cup A^{source})}, g_{iu}, s_i)}\}$. With this, the proof is complete, demonstrating that solving the [MIP] yields an integer optimal solution.